\documentclass[12pt]{elsarticle}
\usepackage{graphicx}
\usepackage{epsfig}
\usepackage{amsmath}
\usepackage{amssymb}
\usepackage{amsfonts}
\usepackage{euscript}
\begin{document}
\bibliographystyle{plain}
\newcommand{\bea}{\begin{eqnarray}}
\newcommand{\eea}{\end{eqnarray}}
\newcommand{\bfm}[1]{{\mbox{\boldmath{$#1$}}}}
\newcommand{\bfmN}{{\mbox{\boldmath{$N$}}}}
\newcommand{\bfmx}{{\mbox{\boldmath{$x$}}}}
\newcommand{\bfmv}{{\mbox{\boldmath{$v$}}}}
\newcommand{\se}{\setcounter{equation}{0}}
\newtheorem{corollary}{Corollary}[section]
\newtheorem{proposition}{Proposition}[section]
\newtheorem{example}{Example}[section]
\newtheorem{definition}{Definition}[section]
\newtheorem{theorem}{Theorem}[section]
\newtheorem{lemma}{Lemma}[section]
\newtheorem{remark}{Remark}[section]
\newtheorem{result}{Result}[section]
\newcommand{\vtwo}{\vskip 4ex}
\newcommand{\vthree}{\vskip 6ex}
\newcommand{\vfour}{\vspace*{8ex}}
\newcommand{\hone}{\mbox{\hspace{1em}}}
\newcommand{\hon}{\mbox{\hspace{1em}}}
\newcommand{\htwo}{\mbox{\hspace{2em}}}
\newcommand{\hthree}{\mbox{\hspace{3em}}}
\newcommand{\hfour}{\mbox{\hspace{4em}}}
\newcommand{\von}{\vskip 1ex}
\newcommand{\vone}{\vskip 2ex}
\newcommand{\n}{\mathfrak{n} }
\newcommand{\m}{\mathfrak{m} }
\newcommand{\q}{\mathfrak{q} }
\newcommand{\aF}{\mathfrak{a} }

\newcommand{\kl}{\mathcal{K}}
\newcommand{\p}{\mathcal{P}}
\newcommand{\Lt}{\mathcal{L}}
\newcommand{\bv}{{\mbox{\boldmath{$v$}}}}
\newcommand{\bc}{{\mbox{\boldmath{$c$}}}}
\newcommand{\bx}{{\mbox{\boldmath{$x$}}}}
\newcommand{\br}{{\mbox{\boldmath{$r$}}}}
\newcommand{\bs}{{\mbox{\boldmath{$s$}}}}
\newcommand{\bb}{{\mbox{\boldmath{$b$}}}}
\newcommand{\ba}{{\mbox{\boldmath{$a$}}}}
\newcommand{\bn}{{\mbox{\boldmath{$n$}}}}
\newcommand{\bp}{{\mbox{\boldmath{$p$}}}}
\newcommand{\by}{{\mbox{\boldmath{$y$}}}}
\newcommand{\bz}{{\mbox{\boldmath{$z$}}}}
\newcommand{\be}{{\mbox{\boldmath{$e$}}}}

\newcommand{\bP}{{\mbox{\boldmath{$P$}}}}

\newcommand{\M}{\mathcal{M}}
\newcommand{\R}{\mathbb{R}}
\newcommand{\Q}{\mathbb{Q}}
\newcommand{\Z}{\mathbb{Z}}
\newcommand{\N}{\mathbb{N}}
\newcommand{\C}{\mathbb{C}}
\newcommand{\xar}{\longrightarrow}
\newcommand{\ov}{\overline}
 \newcommand{\rt}{\rightarrow}
 \newcommand{\om}{\omega}
 \newcommand{\wh}{\widehat }
 \newcommand{\wt}{\widetilde }
 \newcommand{\g}{\Gamma}
 \newcommand{\lm}{\lambda}

\newcommand{\eN}{\EuScript{N}}
\newcommand{\ncom}{\newcommand}
\newcommand{\norm}{\|\;\;\|}
\newcommand{\inp}[2]{\langle{#1},\,{#2} \rangle}
\newcommand{\nrm}[1]{\parallel {#1} \parallel}
\newcommand{\nrms}[1]{\parallel {#1} \parallel^2}
\title{Adaptive Pad\'e-Chebyshev Type Approximation of Piecewise Smooth Functions}
\author{ S. Akansha\footnote{akansha@math.iitb.ac.in} and S. Baskar\footnote{baskar@math.iitb.ac.in}\\
Department of Mathematics, \\Indian Institute of Technology Bombay,\\ Powai, Mumbai - 400076. India.}
\maketitle{}
\begin{center}{\bf Abstract}\end{center}
A piecewise Pad\'e-Chebyshev type (PiPCT) approximation method is proposed to minimize the Gibbs phenomenon in approximating piecewise smooth functions. A theorem on $L^1$-error estimate is proved for sufficiently smooth functions using a decay property of the Chebyshev coefficients. Numerical experiments are performed to show that the PiPCT method accurately captures isolated singularities of a function without using the positions and the types of singularities. Further, an adaptive partition approach to the PiPCT method is developed (referred to as the APiPCT method) to achieve the required accuracy with a lesser computational cost. Numerical experiments are performed to show some advantages of using the PiPCT and APiPCT methods compared to some well-known methods in the literature.

\section{Introduction}
Often in applications, we come across the problem of approximating a non-smooth function, for instance, while approximating the entropy solution of a hyperbolic conservation law that involves shocks and rarefactions, and while compressing an image that involves edges. The main challenge in such an approximation problem is the occurrence of the well known Gibbs phenomenon in the approximant near a singularity of the target function. The methods that are proven to have a higher order of accuracy for smooth target functions tend to have a lower order of convergence near jump discontinuities (Gottlieb and Shu \cite{got-shu_97a}).

Often nonlinear approximation methods are preferred to linear approximation procedures in order to achieve convergence in the vicinity of singularities with a higher order of accuracy (see DeVore \cite{dev_98a}). Multiresolution schemes have been developed based on the Harten's ENO procedure, where one chooses an appropriate set of nodes locally to minimize the effect of the singularity on the approximation (Arandiga {\it et al.} \cite{ara-etal_05a}). Another approach is to construct the approximation using a standard linear approximation method and then use filters or mollifiers to reduce the Gibbs phenomenon (see Tadmor \cite{tad_07a}). By knowing some more information about the singularities of the target function, say, the location and (or) the nature of the singularity, one can develop efficient and accurate methods (see, for instance, Gottlieb and Shu \cite{got-shu_95a}, Driscoll and Fornberg \cite{dri-for_01a} and Lipman and Levin \cite{lip-lev_10a}).

Another nonlinear approach for approximating a piecewise smooth function is a rational approximation, in particular, the Pad\'e approximation. The main advantage of rationalizing a truncated series is to achieve faster convergence in the case of approximating analytic functions (see Baker {\it et al.} \cite{bak_65a, bak-etal_96a}) when compared to polynomial approximations. The advantage of rational approximation also lies in its property of error autocorrection (see Litvinov \cite{lit_03a}), in which the errors in all the intermediate steps compensate each other at the final step of rationalization. The error autocorrection is desirable for approximating piecewise smooth functions. Geer \cite{gee_95a} used a truncated Fourier series expansion of a periodic even or odd piecewise smooth function to construct a Pad\'e approximant, which has been further developed for general functions by Min {\it et al.} \cite{min-etal_07a}. Hesthaven {\it et al.} \cite{hes-etal_06a} proposed a Pad\'e-Legendre interpolation method to approximate a piecewise smooth function. Kaber and Maday \cite{kab-mad_05a} studied the convergence rate of a sequence of Pad\'e-Chebyshev approximation of the sign function. These nonlinear methods minimize the Gibbs phenomenon significantly. However, the convergence at the jump discontinuity is rather slow. In order to accelerate the convergence rate further, Driscoll and Fornberg \cite{dri-for_01a} (also see Tampos {\it et al.} \cite{tam-etal_12a}) developed a singular Pad\'e approximation to a piecewise smooth function using finite series expansion coefficients or function values. Assuming that the jump locations are known, this method captures the discontinuities of a piecewise smooth function very accurately. The methods mentioned above are global approximations. However, in many applications (for instance, constructing numerical schemes for partial differential equations), it is desirable to have a local approximation technique that minimizes the Gibbs phenomenon without the explicit knowledge of the type, magnitude, and location of the singularities of the target function.

In this article, we intend to develop a local approximation technique based on Pad\'e-Chebyshev type approximation that is suitable for approximating a piecewise smooth function. To this end, we propose a piecewise implementation of the Pad\'e-Chebyshev type approximation (denoted by PiPCT approximation) and show numerically that the proposed method captures the isolated singularities (including jump discontinuities) of a non-smooth function accurately with negligible oscillations. 
To approximate a given function $f\in L^1([a,b])$, we consider two vectors $\bfm{n}_p, \bfm{n}_q \in \mathbb{N}^N$, and a partition $P_N$ of the interval $[a,b]$ consisting of $N$ subintervals $I_j$, $j=0,1,\ldots, N-1$. The PiPCT method computes the Pad\'e-Chebyshev type approximation of order $[n_p^j/n_q^j]$ in each subinterval $I_j$, where $n_p^j$ and $n_q^j$ are the components of $\bfm{n}_p$ and $\bfm{n}_q$, respectively. We prove an $L^1$-error estimate of the PiPCT approximation using a decay estimate of the Chebyshev coefficients (see Trefethen \cite{tre_08a}, Xiang {\it et al.} \cite{xia-etal_10a}, and Majidian \cite{maj_17a}). Further, we demonstrate (numerically) the convergence of the PiPCT approximation in the vicinity of jump discontinuity with a higher order of convergence. The key advantage of the PiPCT approximation is that it reduces the Gibbs phenomenon significantly without any knowledge of the location of the jump discontinuity. We compare the results obtained from the PiPCT approximant with the singular Pad\'e approximant \cite{dri-for_01a} and robust Pad\'e-Chebyshev approximant \cite{gon-etal_11a, gon-etal_13a} (both computed on the whole interval, referred to as a global approximation), and find that the PiPCT approximation captures the singularities accurately, which is comparable with the accuracy of the singular Pad\'e approximation.

Our numerical experiments show that the proposed PiPCT method captures the jump discontinuities  accurately. Also, our theoretical study shows that we do not need a piecewise implementation of PCT approximation in a region where the target function is smooth. Hence, we look for an adaptive approach that gives the optimal partition required to approximate the target function with an equivalent accuracy achieved using the PiPCT method on a uniform partition. Such an adaptive algorithm needs the location of the singularities of the target function. Many works have been devoted in the literature to develop singularity indicator, see for instance, Banerjee and Geer \cite{ban-gee_98a}, Barkhudaryan {\it et al.} \cite{bar-etal_07a}, Eckhoff \cite{eck_95a}, Kvernadze \cite{kve_04a}. Since genuine poles appear near singularities (Baker {\it et al.} \cite{bak-etal_61a}), we identify the intervals where the denominator polynomial of the PCT approximation almost vanishes and define these intervals as {\it bad-cells}. We further refine the partition in an isotropic greedy manner (for details, see chapter 3 of \cite{dev-kun_09a}), which results in an adaptive piecewise Pad\'e-Chebyshev type (APiPCT) method. 

In section \ref{Cheb.Approx.sec}, we discuss the decay property of the Chebyshev coefficients of a function. The PiPCT method is proposed  in section \ref{PCT.approximation.sec}, wherein we have also deduced an $L^1$-order of accuracy of the PiPCT approximant from its corresponding  $L^1$-error estimate (obtained in section \ref{Cheb.Approx.sec}). In section \ref{numerical.comparison.sec},  we make a comparative study of the numerical results obtained from the PiPCT method with those obtained from the singular Padé and the robust Padé-Chebyshev methods.  Section \ref{APiPCT.sec}  is devoted to develop the APiPCT method. In section  \ref{numerical.experiment.APiPCT.sec}, we discuss the computational efficiency of the APiPCT approach, where we also implement  the robust Padé-Chebyshev method on the adaptive partition (denoted by the APiRPCT method) and compare the results obtained from the APiRPCT and the APiPCT methods.
Finally, in section \ref{comment.FD.sec}, we comment on the Froissart doublets arising in the PiPCT approximation and their role in the accuracy of the approximant, numerically.
\section{Chebyshev Approximation}\label{Cheb.Approx.sec}
For a given function $f \in L_\omega^2[-1,1],$ with $\omega(x)=1/\sqrt{1-x^2},$ the Chebyshev series representation of $f$ is given by
\begin{equation}\label{CMap.eq}
f(x)\sim \sideset{}{'}\sum_{k=0}^{\infty} c_k T_k(x), ~~ x\in [-1,1],
\end{equation}
where the prime in the summation indicates that the first term is halved, $T_k(x)$ denotes the Chebyshev polynomial of the first kind with degree $k,$ and $c_k$, $k=0,1,\ldots$, are the Chebyshev coefficients given by 
\bea\label{CSRepCoeff.eq}
c_k = \frac{2}{\pi}\langle f, T_k\rangle_w.
\eea
Here, $\langle\cdot,\cdot\rangle_ \omega$ denotes the weighted $L^2$-inner product.

The Chebyshev coefficients \eqref{CSRepCoeff.eq} are approximated using the formula (see \cite{mas-han_03a,riv_74a})
\begin{equation}\label{eq:approxcoeff}
c_{k,n}:=\frac{2}{n}\sum_{l=1}^{n}f(t_l)T_k(t_l) ,~~k=0,1,\ldots,
\end{equation}
where the quadrature points $t_l$, $l=1,2,\ldots, n$, are the Chebyshev points given by
\bea\label{ChebPts.eq}
t_l	= \cos\left(\dfrac{\left(l-0.5\right) \pi}{n}\right), \hspace{1cm} l=1, 2, \ldots,n. 
\eea
We use the notation
$$ \mathsf{C}_{\infty,n}[f](x):=\sideset{}{'}\sum_{k=0}^{\infty} c_{k,n} T_k(x)$$
to denote the Chebyshev series with approximated coefficients, and the truncated series is denoted by
$$ \mathsf{C}_{d,n}[f](x):=\sideset{}{'}\sum_{k=0}^{d} c_{k,n} T_k(x), ~ d=0,1,2,\ldots.$$

To obtain the Chebyshev series representation of a function $f \in L_\omega^2[a,b],$ we use the change of variable 
\begin{equation}\label{eq:bijection1}
	\mathsf{G}(y) = a+(b -a)\frac{(y+1)}{2}, ~~ y\in [-1,1].
\end{equation}
Using the following decay estimate of the Chebyshev coefficients, we can obtain an $L^1$-error estimate (Theorem \ref{thm:errestimate}) for a truncated Chebyshev series approximation of $f$.
\begin{theorem}[Trefethen \cite{tre_13a}, Majidian \cite{maj_17a}]\label{thm:coeffbound}
For some integer $k\ge 0$, let $f$, $f'$, $\ldots,$ $f^{(k-1)}$ be absolutely continuous on the interval $[a,b]$. 
If $V_k:=\|f^{(k)}\|_T < \infty,$ where
\begin{equation}\label{eq:chebnormscale}
\|f\|_T := \int_{0}^{\pi} \left|{f'\left(\mathsf{G}(\cos \theta )\right)}\right|d\theta,
\end{equation}
then for $j\ge k+1$ and for some $s\ge 0$, we have
\begin{equation}\label{eq:coeffbound1}
|c_j| \leq 
\left\{
\begin{array}{@{}l@{\thinspace}l}
\left(\dfrac{b-a}{2}\right)^{2s+1}  \dfrac{2V_k}{\pi\displaystyle \prod_{i = -s}^{s}(j + 2i)},&
\mbox{ if }~ k = 2s,\medskip\\
 \left(\dfrac{b-a}{2}\right)^{2s+2}  \dfrac{2V_k}{\pi\displaystyle \prod_{i = -s}^{s+1}(j + 2i-1)},&
 \mbox{ if }~ k = 2s+1.
 \end{array}
\right. 
\end{equation}
\end{theorem}
\begin{remark}
Note that if $f^{(k)}$ is absolutely continuous, then $V_k$ is precisely the total variation of $f^{(k)}$ and hence in this case, the assumption that $V_k$ is finite implies $f^{(k)}$ is of BV on $[a,b]$.  On the other hand, if $f^{(k)}$ involves a jump discontinuity, then one has to necessarily use distribution derivative of $f^{(k)}$ in computing $V_k$.  Recently, Xiang \cite{xia_18a} obtained a decay estimate based on the total variation of  $f^{(k)}$.
\end{remark}
For the well-known decay estimate of the Chebyshev coefficients of a real analytic function we refer to Rivlin \cite{riv_74a} (also see Xiang {\it et al.} \cite{xia-etal_10a}). The required error estimate can now be proved using the decay estimate \eqref{eq:coeffbound1} of the Chebyshev coefficients:
\begin{theorem}\label{thm:errestimate}
Assume the hypotheses of Theorem \ref{thm:coeffbound}. Then for any given integers $n$ and $d$ such that $n-1\ge k\ge 1$ and $k\le d \le 2n-k-1$, we have
$$\|f - \mathsf{C}_{d,n}[f]\|_1
 \leq \mathbb{T}_{d,n}$$
where,
\begin{enumerate}
\item if $d=n-l$, for some $l =1,2,\ldots, n-k$, then we have
\begin{equation}\label{eq:totalL1err}
\mathbb{T}_{d,n}:= \left\{
\begin{array}{@{}l@{\thinspace}l}
\left(\dfrac{b-a}{2}\right)^{k+2} \dfrac{4V_k}{k \pi } \left(\Pi_{1,1}(n-l) + \Pi_{1,2}(n-l)\right),  & \text{ if } k = 2s,\medskip\\		
\left(\dfrac{b-a}{2}\right)^{k+2} \dfrac{4V_k}{k\pi} \left(\Pi_{0,0}(n-l) + \Pi_{0,1}(n-l)\right), & \text{ if } k = 2s+1, \\
\end{array}
\right.
\end{equation}
\item else if $d = n+l$, for some $l =0,1,\ldots, n-k-1$, then we have
\begin{equation}
\mathbb{T}_{d,n} := \left\{
\begin{array}{@{}l@{\thinspace}l}
\left(\dfrac{b-a}{2}\right)^{k+2} \dfrac{6V_k}{\pi k } \left(\Pi_{1,0}(n-l) + \Pi_{1,1}(n-l)\right),  & \text{ if } k = 2s,\medskip\\		
\left(\dfrac{b-a}{2}\right)^{k+2} \dfrac{6V_k}{k\pi} \left(\Pi_{0,-1}(n-l) + \Pi_{0,0}(n-l) \right), & \text{ if } k = 2s+1, \\
\end{array}
\right.
\end{equation}
for some integer $s\ge 0$, where
\begin{equation}\label{PI.def}
\Pi_{\alpha,\beta}(\eta) := \dfrac{1}{\displaystyle \prod_{i = -s}^{s-\alpha}(\eta+2i+\beta)},~~\alpha=0,1,~~\beta=-1,0,1,2.
\end{equation}
\end{enumerate}
\end{theorem}
\textbf{Proof:} 
We have 
$$
\big\|f-C_{d,n}[f]\big\|_1
\le
\big\|f-C_{d}[f]\big\|_1 + \big\|C_{d}[f]-C_{d,n}[f]\big\|_1.
$$
For estimating the second term on the right hand side of the above inequality, we use  the well-known result
(see Fox and Parker \cite{fox-par_68a}) 
$$c_{d,n} - c_d = \sum_{j=1}^\infty (-1)^j\big(c_{2jn-d} + c_{2jn+d}\big),$$
for $0\le d< 2n$.  Using this property, with an obvious rearrangements of the terms in the series, we can obtain
\begin{align}\label{Est01.eq}
\|f - \mathsf{C}_{d,n}[f]\|_1
& \leq (b-a)\mathcal{E},
\end{align}
where
$$\mathcal{E}:=\left\{\displaystyle{\sum_{j=d+1}^\infty} |c_j|+\sum_{j=1}^\infty \sum_{i=2jn-d}^{2jn+d}|c_i|\right\}.$$
By adding some appropriate positive terms, we can see that (also see Xiang {\it et al.} \cite{xia-etal_10a})
\begin{equation}\label{E.estimate}
\mathcal{E} \le 
\left\{\begin{array}{ll}
2\displaystyle{\sum_{i=n-l+1}^{\infty}}|c_i|,&\mbox{for } d=n-l,~l =1,2,\ldots, n,\\
3\displaystyle{\sum_{i=n-l}^{\infty}}|c_i|,&\mbox{for } d=n+l,~l=0,1,\ldots,n-1.
\end{array}\right.
\end{equation}
We restrict the integer $d$ to $k\le d \le 2n-k-1$ so that the decay estimate in Theorem \ref{thm:coeffbound} can be used.   Now using the telescopic property of the resulting series (see also Majidian \cite{maj_17a}), we can arrive at the required estimate.
\qed
\begin{remark}\label{UB.behavior.remark}
From the above theorem, we see that for a fixed $n$ (as in the hypotheses), the upper bound $\mathbb{T}_{d,n}$ decreases for $d<n$ and increases for $d\ge n$.  Further, we see that $\mathbb{T}_{n-l-1,n} = \frac{2}{3}\mathbb{T}_{n+l,n}$,
however computationally  $\mathsf{C}_{n-l-1,n}[f]$ is more efficient than $\mathsf{C}_{n+l,n}[f].$ A similar situation occurs in the case of  Pad\'e-Chebyshev type approximation as shown in subsection \ref{higherPCT.lowerPCT.ssec} and numerically shown in  \ref{degree.adaptation.subs}.\qed
\end{remark}
On similar lines, one can prove the error estimate of a real analytic function $f$. See Xiang {\it et al.} \cite{xia-etal_10a} for more details.  

Note that the hypotheses of the above theorems enable us to use the error estimate only for functions that are at least continuous. It is well-known that if a function has a jump discontinuity, then the Chebyshev approximant may develop Gibbs phenomenon in the vicinity of the jump. One can reduce the Gibbs phenomenon significantly (if not entirely) using a rational approximation, for instance, the Pad\'e-Chebyshev approximation.  
\section{Pad\'e-Chebyshev type approximation}\label{PCT.approximation.sec}
In this section, we recall the basic construction of the Pad\'e-Chebyshev type approximation of a given function $f\in L^2_w[-1,1]$ and further propose a piecewise implementation of this approximation (in subsection \ref{piecewise.implementation.ssec}).

For $x\in [-1,1],$ we use the notation 
\bea\label{CComplexNotation.eq}
\mathsf{C}_{\infty}[f](z)
:= \sideset{}{'}\sum_{k=0}^{\infty} c_k z^k;~~z = e^{i\cos^{-1}(x)},~ i=\sqrt{-1},
\eea
which is a power series defined on the unit circle centered at the origin in the complex plane with real coefficients given by \eqref{CSRepCoeff.eq}.  We can see that the real part of this series is precisely the Chebyshev series given by \eqref{CMap.eq}.  

 For any given integers $n_p\ge n_q \geq 1$, a rational function 
 \bea\label{PadeApprox.eq}
 \mathsf{R}_{n_p,n_q}(z) := \dfrac{P_{n_p}(z)}{Q_{n_q}(z)}
 \eea
  with numerator polynomial $P_{n_p}(z)$ of degree  $\leq n_p$  and denominator polynomial $Q_{n_q}(z)$ of degree $\leq n_q$ with $Q_{n_q} \neq 0$ satisfying (see \cite{bec-mat_15a,dri-for_01a,tam-etal_12a})
\begin{equation} \label{eq:P_linear}
Q_{n_q}(z)\mathsf{C}_{\infty}[f](z)-P_{n_p}(z)= \mathit{O}(z^{n_p+n_q+1}), \hspace*{.5cm} z\rightarrow 0,
\end{equation}
is a Pad\'e approximant of $\mathsf{C}_{\infty}[f](z)$ of order $[n_p/n_q]$. Such a Pad\'e approximation exists, and the real part of $ \mathsf{R}_{n_p,n_q}(z)$ is an approximation of $f(x)$, which is referred to as a {\it Pad\'e-Chebyshev type} (PCT) approximant of $f$ (see \cite{tam-etal_12a}).

From \eqref{eq:P_linear}, we see that the coefficients of the denominator polynomial $Q_{n_q}(z)$ are obtained as a solution of the $n_q\times (n_q+1)$ homogeneous Toeplitz system \cite{dri-for_01a,hes-etal_06a,tam-etal_12a,bec-mat_15a}
\begin{equation*}
\begin{bmatrix}
c_{n_p+1} & c_{n_p} &\cdots & c_{n_p-n_q+1} \\
c_{n_p+2} & c_{n_p+1} &\cdots & c_{n_p-n_q+2} \\
\vdots &\vdots&\vdots &\vdots \\
c_{n_p+n_q} &c_{n_p+n_q-1} & \cdots & c_{n_p} 
\end{bmatrix} \begin{bmatrix}
q_0 \\ q_1 \\ \vdots \\ q_{n_q}
\end{bmatrix}
= \begin{bmatrix}
0 \\ 0 \\ \vdots \\ 0
\end{bmatrix},
\end{equation*}
which has a non-trivial solution.

We write the above system as
\begin{equation}\label{eq:linearsystem}
A_{n_p,n_q}\bfm{q} = \bfm{0},
\end{equation}
where $A_{n_p,n_q}=(a_{r,s})$ is the Toeplitz matrix with
$$a_{r,s}=c_{n_p+r-s+1},~~r=1,2\ldots,n_q, s=1,\ldots,n_q+1, $$
and $\bfm{q} = (q_0,q_1,\ldots,q_{n_q})^T.$

Once the coefficients of the denominator polynomial are known, we can compute the coefficients of the numerator polynomial  $P_{n_p}(z)$ using the following matrix vector multiplication:
\begin{equation}
\begin{bmatrix}
p_0 \\ p_1 \\ \vdots\\ p_{n_q}\\ \vdots \\ p_{n_p}
\end{bmatrix} =
\begin{bmatrix}
c_{0}/2 & 0 & \cdots & 0 \\
c_{1} &c_0/2 & \cdots & 0 \\
\vdots & & & \\
c_{n_q} & c_{n_q-1} & \cdots & c_{0}/2 \\
\vdots & & & \\
c_{n_p} & c_{n_p-1} & \cdots & c_{n_p-n_q} 
\end{bmatrix} \begin{bmatrix}
q_0 \\ q_1 \\ \vdots \\ q_{n_q}
\end{bmatrix} \label{eq:PCnum}.
\end{equation}
If $A_{n_p,n_q}$ is of full rank, then a normalization $q_0 = 1$ leads to a unique (normalized) PCT approximant of $f$ of order $[n_p/n_q]$. 

A Pad\'e-Chebyshev type approximation of $f$ of order $[n_p/n_q]$ can be computed for a given set of Chebyshev coefficients $\{c_0,c_1,\ldots,c_{n_p+n_q}\}$.  Since these coefficients often cannot be obtained exactly, we use the approximated coefficients $\{c_{0,n},c_{1,n},$
$\ldots,c_{n_p+n_q,n}\}$ given by the formula \eqref{eq:approxcoeff} to compute the PCT approximant and denote it by $ \mathsf{R}_{n_p,n_q}^{n}.$
We denote the corresponding Toeplitz matrix in \eqref{eq:linearsystem} by  $A_{n_p ,n_q}^n$. 

\subsection{Error Estimate}
Our aim in this subsection is to use the decay estimate \eqref{eq:coeffbound1} and obtain an error estimate for the Pad\'e-Chebyshev type approximation of the truncated Chebyshev series representation of a function $f$ with approximate coefficients obtained by using the Chebyshev quadrature rule with $n$ quadrature points. 

We truncate the Chebyshev series up to $d=n_p+n_q$ number of terms and take the approximate Chebyshev coefficients in \eqref{eq:P_linear} to write 
$$
Q_{n_q}(z)\mathsf{C}_{d,n}[f](z)-P_{n_p}(z)= 
a_{d+1} z^{d+1} + \ldots + a_{d+n_q} z^{d+n_q},
$$
where
$$
a_{d+j} = \sum_{i=j}^{n_q} q_i c_{d-i+j,n}, ~~j=1,2,\ldots,n_q.
$$
Dividing both sides by $Q_{n_q}(z)$, taking modulus on both sides, and noting that $|z|\le 1$, we can write
\bea\label{Error.Estimate.eq}
\left|\mathsf{C}_{d,n}[f](z)-\mathsf{R}_{n_p,n_q}^{n}(z)\right|
			\le
	\dfrac{n_q\displaystyle{q^*}}{\displaystyle{\prod_{j=1}^{n_q}\big|(z-z_j)\big|}}
		\sum_{i=1}^{n_q} |c_{n_p+i,n}|
\eea
where $q^* = \left|\frac{1}{q_{n_q}}\right|\displaystyle{\max_{1\le i\le n_q}} |q_i|$ and $z_j$, $j=1,2,\ldots,n_q$, are the zeros of the polynomial $Q_{n_q}(z)$.

Again since
$$c_{n_p+i,n}
 	= c_{n_p+i} + \sum_{j=1}^\infty (-1)^j\big(c_{2jn-n_p-i} + c_{2jn+n_p+i}\big),$$
we can write
\bea\label{ApproxCoeff.Estimate.eq}
 \sum_{i=1}^{n_q}|c_{n_p+i,n}|
\le  \sum_{j=n_p+1}^{d}|c_{j}| 
		+ \sum_{j=1}^\infty 
		\sum_{i=2jn-d}^{2jn+d}\big|c_{i}\big|.
\eea
Let us now consider two cases, namely, $0\le d<n$ and $n\le d<2n$.  In the first case, let us take $d=n-l$ for $l=1,2,\ldots n.$ By adding some appropriate positive terms on the right hand side of \eqref{ApproxCoeff.Estimate.eq}, we can write
$$\sum_{i=1}^{n_q}|c_{n_p+i,n}| \le \sum_{j=n_p+1}^{\infty}|c_{j}|,~~\mbox{for}~~ 0\le d<n. $$
For the second case, let us take $d=n+l$ for $l=0,1,\ldots n-1.$ Then the second term on the right hand side of \eqref{ApproxCoeff.Estimate.eq} can be written as
$$\sum_{j=1}^\infty \sum_{i=2jn-d}^{2jn+d}|c_i|
=\left(\sum_{i=n-l}^{n+l}|c_i| + \sum_{i=3n-l}^{3n+l}|c_i| + \ldots \right) + 
\sum_{i=n+l+1}^{\infty}|c_i| .
$$
Hence, \eqref{ApproxCoeff.Estimate.eq} can be written as
$$
 \sum_{i=1}^{n_q}|c_{n_p+i,n}|
\le
 \sum_{j=n_p+1}^{n+l}|c_{j}| 
		+\left(\sum_{i=n-l}^{n+l}|c_i| + \sum_{i=3n-l}^{3n+l}|c_i| + \ldots \right) + 
\sum_{i=n+l+1}^{\infty}|c_i| 
$$
Let $n^* = \min\{n_p+1,n-l\}$.  Then, adding some appropriate positive terms, the above inequality can be written as
$$
 \sum_{i=1}^{n_q}|c_{n_p+i,n}|
\le 2\sum_{i=n^*}^{\infty}|c_i|.
$$
Combining the above two inequalities, we can write
\begin{equation}\label{CC.estimate}
 \sum_{i=1}^{n_q}|c_{n_p+i,n}|
\le  \left\{\begin{array}{ll}
\displaystyle{\sum_{i=n_p+1}^{\infty}}|c_i|,&\mbox{for } d=n-l,~l=1,2,\ldots,n\\
2\displaystyle{\sum_{i=n^*}^{\infty}}|c_i|,&\mbox{for } d=n+l,~l=0,1,\ldots,n-1,
\end{array}\right.
\end{equation}
By substituting \eqref{CC.estimate}, we get the right hand side of \eqref{Error.Estimate.eq} involving the exact Chebyshev coefficients, and hence the decay estimate \eqref{eq:coeffbound1} in Theorem \ref{thm:coeffbound} can be used, provided all the hypotheses of the theorem are satisfied.
\begin{theorem}
\label{PCEstimate.thm}
Assume the hypotheses of Theorem \ref{thm:coeffbound} with $k>1$. Choose the integers $n_0\ge k+2$,$n_p\ge k+1$, $n_q\ge 1$, and $d=n_p+n_q$. Then for a given $n\ge n_0$ and $k+2\le d\le n+n_0-k-1,$ we can find a set $\Gamma_{n_q}\subset \mathbb{C}$ such that
\begin{equation}\label{ErrEstFinal.eq}
\left|\mathsf{C}_{d,n}[f](z)-\mathsf{R}_{n_p,n_q}^{n}(z)\right|
< \frac{1}{(2n_q)^\nu}\left(\dfrac{4V_k\displaystyle{ q^*}}{k\pi}\right)\left(\dfrac{b-a}{2}\right)^{k+1-\nu},
\end{equation}
for any $0<\nu < 1$ and for all $z\in \Gamma_{n_q}^c \cap \{z\le 1\}$.
\end{theorem}
{\bf Proof:}
Let us first consider the case when $k+2\le d<n$.  Since $n_p\ge k+1$, we can use the decay estimate  in Theorem \ref{thm:coeffbound} to get
$$
\displaystyle{\sum_{i=n_p+1}^{\infty}}|c_i|
\le
\left\{
\begin{array}{@{}l@{\thinspace}l}
\left(\dfrac{b-a}{2}\right)^{k+1} \dfrac{V_k}{k\pi} \Big\{\Pi_{1,1}(n_p)
+\Pi_{1,2}(n_p)
\Big\},&
\mbox{ if }~ k = 2s,\medskip\\
 \left(\dfrac{b-a}{2}\right)^{k+1} \dfrac{V_k}{k\pi} \Big\{\Pi_{0,0}(n_p)
+\Pi_{0,1}(n_p)\Big\},&
 \mbox{ if }~ k = 2s+1.
 \end{array}
\right. 
$$
Let us now consider  $n\le d \le n+n_0-k-1$. Since $n_p\ge k+1$ and  $d\le n+n_0-k-1,$ we see that $n^*\ge k+1$.  Hence, we can use Theorem \ref{thm:coeffbound} to get
\begin{equation*}
\displaystyle{\sum_{i=n^*}^{\infty}}|c_i|
\le
{\small \left\{
\begin{array}{@{}l@{\thinspace}l}
\left(\dfrac{b-a}{2}\right)^{k+1}\dfrac{V_k}{k\pi} \Big\{\Pi_{1,0}(n^*)
+\Pi_{1,1}(n^*)
\Big\},&
\mbox{ if }~ k = 2s,\medskip\\
 \left(\dfrac{b-a}{2}\right)^{k+1}\dfrac{V_k}{k\pi}\Big\{
					\Pi_{0,-1}(n^*)
					+\Pi_{0,0}(n^*)
				\Big\},&
 \mbox{ if }~ k = 2s+1.
 \end{array}\right.}
\end{equation*}
Further, we can obtain the estimate
\begin{equation}\label{PIE.def}
\Pi_{\alpha,\beta}(\eta) \le \dfrac{1}{(\eta-k)^k},~~
\end{equation}
where 
$$
\eta = \left\{\begin{array}{ll}
			n_p,&\mbox{ for } k+2\le d<n,\\
			n^*,&\mbox{ for } n\le d\le n+n_0-k-1,
	\end{array}\right.
$$
$\alpha=0,1,$ $\beta=-1,0,1,2,$ and $n^* = \min\{n_p+1,n-l\}$, $l=0,1,\ldots,n_0-k-1.$

Using the above  estimates in \eqref{CC.estimate} and substituting the resulting estimate in \eqref{Error.Estimate.eq}, we get
\begin{equation}\label{Error.Estimate.V3.Simple.eq}
\left|\mathsf{C}_{d,n}[f](z)-\mathsf{R}_{n_p,n_q}^{n}(z)\right|\le \dfrac{n_q\displaystyle{q^*}}{\displaystyle{\prod_{i=1}^{n_q}\big|(z-z_i)\big|}} \left(\dfrac{b-a}{2}\right)^{k+1} \dfrac{V_k}{k\pi} 
		 \dfrac{4}{(\eta-k)^k}.
\end{equation}

All remains now is to estimate the denominator polynomial on the right hand side of  the above inequality. For this, we use the well-known Boutroux-Cartan lemma (see Boas \cite{boa_54a} and Wallin \cite{wal_74a}). 

For a fixed real number $0<\nu<1$, let us define
\begin{equation}\label{Delta.exp}
\delta^{n_q} = \frac{n_q^{\nu+1}(b-a)^\nu}{(\eta-k)^k}>0.
\end{equation}
By the Boutroux-Cartan lemma, there exists at most $n_q$ disks $C_{r_j}$, $j=1,2,\ldots,m (\le n_q)$, each with  radius $r_j$ satisfying
$\displaystyle{\sum_{j=1}^m} r_j \le 2e\delta,$
such that
$$\left(\prod_{i=1}^{n_q}\big|(z-z_i)\big|\right)^{-1} < \left(\frac{1}{\delta}\right)^{n_q},$$
for all $z\in \Gamma_\delta^c$, where 
$z_i\in \Gamma_\delta := \bigcup_{j=1}^m C_{r_j}, ~~i=1,2,\ldots,n_q.$ The inequality \eqref{ErrEstFinal.eq} follows by substituting the above estimate in \eqref{Error.Estimate.V3.Simple.eq}.
\begin{remark}
We can observe that the right hand side of \eqref{ErrEstFinal.eq} tends to zero as $n_q\rightarrow \infty$.  Further, since $\delta^{n_q}\rightarrow 0$ as $n_q\rightarrow \infty$, the Boutroux-Cartan exceptional set $\Gamma_\delta$ tends to a countable set. 
\end{remark}
\subsection{Computing higher order PCT approximants using lower order PCT approximants}\label{higherPCT.lowerPCT.ssec}
This subsection presents an interesting property of $\mathsf{R}_{n_p,n_q}^{n}$ (also see Cuyt and Wuytack \cite{cuy-wuy_87a} for a similar result) which motivated us in choosing a suitable degree of the numerator polynomial in subsection \ref{degree.adaptation.subs}.
\begin{proposition}\label{prop:1}
	Let $f\in L^{2}_{\omega}[-1,1]$ and consider the approximated Chebyshev series $f(x) \approx \sum_{k=0}^{\infty}c_{k,n}T_k(x)$ for $x\in [-1,1].$ Let $\mathsf{R}_{n_p,n_q}^{n}(z)$ be the unique Pad\'e approximation of $\mathsf{C}_{\infty,n}[f](z)$ of order $[n_p/n_q]$. Then
for Pad\'e approximations of order $[(n-j-1)/n_q]$ and $[(n+j)/n_q]$, $j = 0,1,\ldots,n-n_q-1$, $n_q<n-1$, the corresponding denominator coefficient vectors $\bfm{q}_{n_q}^{n-j-1}$ and $\bfm{q}_{n_q}^{n+j}$ satisfy
		\begin{equation}\label{eq:Denrelationgen}
		\bfm{q}_{n_q}^{n-j-1} = R_{n_q+1}\bfm{q}_{n_q}^{n+j},
		\end{equation}
where $R_{n_q+1}$ is the `flip' matrix of size $(n_q+1) \times (n_q+1)$ with $1$ on anti-diagonal and $0$ elsewhere.
\end{proposition}
We use the following anti-symmetric property of the approximated Chebyshev series coefficients to prove Proposition \ref{prop:1}.
\begin{lemma}\label{lem:2}
	The approximated Chebyshev series coefficients \eqref{eq:approxcoeff} of $f$ using $n$ quadrature points satisfy, for $k=1,3,5,\ldots$,
	\begin{equation}\label{eq:result1}
	\left.\begin{array}{rlr}
	c_{kn+j,n}  &= -c_{|kn-j|,n}, & \hspace{.3cm} j = 1,2,\ldots,2n, \\
	c_{kn,n} &= 0.&
	\end{array}\right\}
	\end{equation}
\end{lemma}
\textbf{Proof: } The result follows immediately from the identity
$$T_{kn}(t_l) T_j(t_l)= \frac{1}{2}\Big(T_{kn+j}(t_l) + T_{|kn-j|}(t_l)\Big),~~l=1,2,\ldots,n.$$\qed

\noindent\textbf{Proof of Proposition \ref{prop:1}: }
By the construction of the Pad\'e-Chebyshev approximants of order $[(n-j-1)/n_q]$ and $[n+j/n_q]$, for $j = 0,1,\ldots,n-n_q-1$, the denominator coefficient vectors $\bfm{q}_{n_q}^{n-j-1}$ and $\bfm{q}_{n_q}^{n+j}$ are obtained as solutions of the linear systems
$$
A_{n-j-1,n_q}^n\bfm{q}_{n_q}^{n-j-1} = \bfm{0}
$$
and
$$
A_{n+j,n_q}^n\bfm{q}_{n_q}^{n+j} = \bfm{0},
$$
respectively.  

Using Lemma \ref{lem:2} we can write $A_{n+j,n_q}^n = (b_{r,s})$, where
$$b_{r,s} = -c_{|n-(j+r-s+1)|,n},~~r=1,2\ldots,n_q, s=1,\ldots,n_q+1.$$
This shows that  by flipping rows and columns of the matrix $A_{n-j-1,n_q}^n,$ and then multiplying by $-1$, we obtain the matrix $A_{n+j,n_q}^n$.  That is,
\begin{equation}\label{eq:prop1}
	A_{n+j,n_q}^n = -R_{n_q}A_{n-j-1,n_q}^nR_{n_q+1},
\end{equation}
where $R_{n_q}$ and $R_{n_q+1}$ are `flip' matrices of size $n_q$ and $n_q+1$, respectively.  

Let $\bfm{x} \in \ker (A_{n+j,n_q}^n)$. Then from \eqref{eq:prop1}, we see that 
$$A_{n-j-1,n_q}^nR_{n_q+1}\bfm{x}  = \bfm{0}.$$	
Therefore, for any vector $\bfm{x} \in \ker(A_{n+j,n_q}^n)$, we have $R_{n_q+1}\bfm{x} \in \ker(A_{n-j-1,n_q}^n)$ and vice versa.	This completes the proof.\qed
\subsection{Construction of PiPCT Approximation}\label{piecewise.implementation.ssec}
Let $f\in L^1[a,b]$. Discretize the interval $I:=[a,b]$ into $N$ subintervals (need not be equally spaced),  denoted by $I_j:=[a_j,b_j]$, $j=0,1,\ldots,N-1$, where
$$a=a_0<b_0=a_1<b_1=a_2<\cdots<b_{N-2}=a_{N-1}<b_{N-1}=a_N=b$$
and denote the partition by $P_N:=\{a_0,a_1,\ldots,a_N\}.$

For any given integers $n$ and $N$, and $N$-tuples $\bfm{n_p}=(n_{p}^{0},\ldots,n_{p}^{N-1})$ and $\bfm{n_{q}}=(n_{q}^{0},\ldots,n_{q}^{N-1}),$   
construct the PiPCT approximation of $f$ as follows:\\
\begin{enumerate}
\item Generate $n$ Chebyshev points $\{t_l~:~ l=1,2,\cdots,n\},$ given by \eqref{ChebPts.eq}, in the reference interval $[-1,1]$. \\
\item For each $j=0,1, \ldots, N-1$, consider the bijection map $\mathsf{G}_j : [-1,1] \rightarrow I_j$ given by
	\begin{equation}\label{eq:bijection}
	\mathsf{G}_j(y) = a_j+(b_j -a_j)\frac{(y+1)}{2}.
	\end{equation}
Obtain the approximate Chebyshev coefficients in the $j^{\rm th}$ cell, denoted by $c_{k,n}^j,$ for $k=0,1,\cdots,n_{p}^{j}+n_{q}^{j},$ using the formula \eqref{eq:approxcoeff} with the values of $f$ evaluated at $\mathsf{G}_j(t_l)$, $l=1,2,\cdots,n$.\\
\item For each $j=0,1,\ldots,N-1$, compute the PCT approximant $\mathsf{R}_{n_{p}^j,n_{q}^j}^{n}$ of order $[n_p^j/n_q^j]$ of the function $f|_{I_j}$.\\
\item Define the {\it piecewise Pad\'e-Chebyshev type}  (PiPCT) approximation of $f$ in the interval $I$ with respect to the given partition as
\bea\label{PPC.eq}
\mathsf{R}^{n,N}_{\mbox{\tiny $\bfm{n}$}_p,\mbox{\tiny $\bfm{n}$}_q}(x):=
\begin{cases}
\mathsf{R}_{n_{p}^0,n_{q}^0}^{n}(x),&\text{if }x \in I_0,\medskip\\
\mathsf{R}_{n_{p}^1,n_{q}^1}^{n}(x) ,&\text{if } x\in I_1,\medskip\\
\vdots \\
\mathsf{R}_{n_{p}^{N-1},n_{q}^{N-1}}^{n}(x),&\text{if } x \in I_{N-1}.
\end{cases}
\eea
\end{enumerate}
Using Theorems \ref{thm:errestimate} and \ref{PCEstimate.thm}, we can obtain an error bound for the PiPCT approximation.
\begin{theorem}\label{errorestimate.thm}
Assume that $f:[a,b]\rightarrow \mathbb{R}$ is a piecewise smooth function.  For some integer $k> 1$, let
\begin{enumerate}
\item $f$, $f'$, $\ldots,$ $f^{(k)}$ be absolutely continuous on $[a,b]$; and
\item 
$V_{k}:=\|f^{(k)}\|_T < \infty$.
\end{enumerate}
Let $P_N$ be a partition of $[a,b]$, for some integer $N\ge 1$, and $n_0\ge k+2$ and $n\ge n_0$ be integers.
Then for any choice of  ${n}_{p}^{j}\ge k+1$ and $1\le {n}_{q}^{j}\le{n}_{p}^{j}$, $j=1,2,\ldots, N-1$, such that  $k+2\le d_*\leq d^* \le n+n_0-k-1$, where $d_* := \min\{n_p^j+n_q^j~|~j=0,1,\ldots,N-1\}$ and $d^*:= \max\{n_p^j+n_q^j~|~j=0,1,\ldots,N-1\}$, there exists a constant $K>0$ and a set $\mathcal{I}_{\mbox{\tiny $\bfm{n}$}_q}\subset [a,b]$ such that
\begin{equation}\label{Error.Est.PiPCT}
 \big\|f-\mathsf{R}^{n,N}_{\mbox{\tiny $\bfm{n}$}_p,\mbox{\tiny $\bfm{n}$}_q}\big\|_{1,\mbox{\tiny $\bfm{n}$}_q}
< K(b-a)\Delta \left(\dfrac{h}{2}\right)^{k+1-\nu},
\end{equation}
for any $0<\nu<1$,
where $\|\cdot\|_{1,\mbox{\tiny $\bfm{n}$}_q}$ denotes the $L^1$-norm taken over the set $\mathcal{I}_{\mbox{\tiny $\bfm{n}$}_q}$, 
$$h = \max\Big\{b_j-a_j ~|~ j=0,1,\ldots,N-1\Big\}$$
and 
$\Delta =\max\left\{
\Delta^*,1/(2{n}_{q}^{*})^\nu
\right\}$,
with ${n}_{q}^{*} = \min_j({n}_{q}^{j})$ and
$$\Delta^* =  \max_{j=0,\ldots,N-1}\left(\left\{\begin{array}{ll}
				\dfrac{1}{(d^j-k)^k},&\mbox{ if } k+2\le d^j<n,\medskip\\
				\dfrac{1}{(n-n_0+1)^k},&\mbox{ if } n\le d^j\le n+n_0-k-1.					
		\end{array}\right.\right).
$$
\end{theorem}
\textbf{Proof: } 
For each $j=0,1,2,\ldots,N-1$, let us first find $\mathbb{T}_{d^j,n}$, such that
$$
\sum_{j=0}^{N-1}\big\|f|_{I_j}-\mathsf{C}_{d^j,n}[f|_{I_j}]\big\|_1\le \sum_{j=0}^{N-1}\mathbb{T}_{d^j,n},
$$
where $\|\cdot\|_1$ denotes the $L^1$-norm taken over the interval $[a,b]$ and $d^j = n_{p}^j+n_{q}^j$.

Since $\bfm{n}_p$ and $\bfm{n}_q$ are chosen such that that $k+2\le d_*\leq d^* \le n+n_0-k-1$,  we see that for each $j=0,1,\ldots,N-1$, we have $k\le d^j\le 2n-k-1$.  Thus, we can use  Theorem \ref{thm:errestimate} to get
\begin{enumerate}
\item for $d^j=n-l$, where $ l =1,2,\ldots, n-k-2$, and for some integer $s\ge 1$, we have
\begin{equation}\label{eq:totalL1err}
\mathbb{T}_{d^j,n}= \left\{
\begin{array}{@{}l@{\thinspace}l}
\left(\dfrac{b_j-a_j}{2}\right)^{k+2} \dfrac{4V_{k}^j}{k \pi } \left(\Pi_{1,1}^j(n-l) + \Pi_{1,2}^j(n-l)\right),  & \text{ if } k = 2s,\medskip\\		
\left(\dfrac{b_j-a_j}{2}\right)^{k+2} \dfrac{4V_{k}^j}{k\pi} \left(\Pi_{0,0}^j(n-l) + \Pi_{0,1}^j(n-l)\right), & \text{ if } k = 2s+1, \\
\end{array}
\right.
\end{equation}
\item and for $d^j = n+l,$ where $l =0,1,\ldots, n_0-k-1$, and for some integer $s\ge 1$, we have
\begin{equation}\label{eq:totalL1err2}
\mathbb{T}_{d^j,n} = \left\{
\begin{array}{@{}l@{\thinspace}l}
\left(\dfrac{b_j-a_j}{2}\right)^{k+2} \dfrac{6V_{k}^j}{\pi k } \left(\Pi_{1,0}^j(n-l) + \Pi_{1,1}^j(n-l)\right),  & \text{ if } k = 2s,\medskip\\		
\left(\dfrac{b_j-a_j}{2}\right)^{k+2} \dfrac{6V_{k}^j}{k\pi} \left(\Pi_{0,-1}^j(n-l) + \Pi_{0,0}^j(n-l) \right), & \text{ if } k = 2s+1, \\
\end{array}
\right.
\end{equation}
\end{enumerate}
For each $j=0,1,\ldots, N-1$, we have the following estimate
$$
\Pi_{\alpha,\beta}^j						
\le\left\{\begin{array}{ll}
				\dfrac{1}{(d^j-k)^k},&\mbox{ if } k+2\le d^j<n,\medskip\\
				\dfrac{1}{(n-n_0+1)^k},&\mbox{ if } n\le d^j\le n+n_0-k-1.							\end{array}\right.
$$
Therefore, we can write
\begin{equation}\label{Term1.Estim.eq}
\sum_{j=0}^{N-1}\mathbb{T}_{d^j,n} \le
\left(\dfrac{h}{2}\right)^{k+1} (b-a)K'\Delta^*.
\end{equation}
where 
$$K'=\max_{j}\left( \dfrac{6V_{k}^j}{k\pi}\right).$$

Next, let us use  Theorem \ref{PCEstimate.thm} to obtain an estimate for the error in the PiPCT when compared to the truncated Chebyshev series in each subinterval $I_j$,  for $j=0,1,2,\ldots,N-1$.

For a given sub-interval $I_j$, let us define the set
$$
\mathcal{I}_j = \Big\{t\in I_j~|~\exp\big(i\cos^{-1}\big(\mathsf{G}^{-1}_j(t)\big)\big)\in \Gamma_{n_{q}^j}^c \Big\},
$$
where $\Gamma_{n_{q}^j}\subset \mathbb{C} $ is the Boutroux-Cartan exceptional set as obtained in Theorem \ref{PCEstimate.thm} for the polynomial $Q_{n_{q}^j}(z)$.

We apply Theorem \ref{PCEstimate.thm} in each subinterval to get
$$
\left|\mathsf{C}_{d^j,n}[f|_{I_j}](z)-\mathsf{R}_{n_{p}^j,n_{q}^j}^{n}(z)\right|
<
 \frac{1}{(2n_{q}^{j})^\nu}\left(\dfrac{4V_{k}^{j}\displaystyle{ q_{j}^*}}{k\pi}\right)\left(\dfrac{b_j-a_j}{2}\right)^{k+1-\nu},
$$
for any $0<\nu<1$ and  for all $z\in  \Gamma_{n_{q}^j}^c \cap \{z\le 1\}$.
Therefore, we can write
\begin{equation}
\left|\mathsf{C}_{d^j,n}[f|_{I_j}](t) - \mathsf{R}_{n_{p}^j,n_{q}^j}^{n}(t)\right|
< 
 \frac{1}{(2n_{q}^{j})^\nu}\left(\dfrac{4V_{k}^{j}\displaystyle{ q_{j}^*}}{k\pi}\right)\left(\dfrac{b_j-a_j}{2}\right)^{k+1-\nu},\label{Term2.Estim.eq}
\end{equation}
for all $t\in \mathcal{I}_j$.  Define
$$\mathcal{I}_{\mbox{\tiny $\bfm{n}$}_q} = \bigcup_{j=1}^{N-1}\mathcal{I}_j.$$
We have
\begin{eqnarray}
\big\|f-\mathsf{R}^{n,N}_{\mbox{\tiny $\bfm{n}$}_p,\mbox{\tiny $\bfm{n}$}_q}\big\|_{1,\mbox{\tiny $\bfm{n}$}_q} 
&=&\int_{\mathcal{I}_{\mbox{\tiny $\bfm{n}$}_q}}
|f(t)-\mathsf{R}^{n,N}_{\mbox{\tiny $\bfm{n}$}_p,\mbox{\tiny $\bfm{n}$}_q}(t)| dt \nonumber\\
&\le& \sum_{j=0}^{N-1}\big\|f|_{I_j}-\mathsf{C}_{d^j,n}[f|_{I_j}]\big\|_1 +  \sum_{j=0}^{N-1}\int_{\mathcal{I}_j}\left|\mathsf{C}_{d^j,n}[f|_{I_j}](t) - \mathsf{R}_{n_{p}^j,n_{q}^j}^{n}(t)\right| dt.\nonumber
\end{eqnarray}

Using the estimates \eqref{Term1.Estim.eq} and \eqref{Term2.Estim.eq} in the above inequality, we get
\begin{eqnarray*}
\big\|f-\mathsf{R}^{n,N}_{\mbox{\tiny $\bfm{n}$}_p,\mbox{\tiny $\bfm{n}$}_q}\big\|_{1,\mbox{\tiny $\bfm{n}$}_q}
&<& \left(\dfrac{h}{2}\right)^{k+1} (b-a)K'\Delta^* + 
 \frac{K''}{(2{n}_{q}^{*})^\nu}\left(\dfrac{h}{2}\right)^{k+1-\nu}(b-a),
\end{eqnarray*}
 where
$$K''=\left(\dfrac{4\max_{j}(V_{k}^{j}\displaystyle{ q_{j}^*})}{k\pi}\right) 
\mbox{ and } {n}_{q}^{*} = \min_j({n}_{q}^{j}).$$
Take $K=\max(K',K'')$ and $\Delta = \max(\Delta^*,1/(2{n}_{q}^{*})^\nu)$ to get the desired result.
\qed

We observe the following two properties of PiPCT approximation from the above theorem:
\begin{enumerate}
\item The Boutroux-Cartan exceptional set tends to a countable set as $N\rightarrow \infty$ as well as $\bfm{n}_q\rightarrow \infty$.
\item The estimate \eqref{Error.Est.PiPCT} shows that the $L^1$-rate of convergence of the PiPCT approximation is at least $k$ as $N\rightarrow \infty$. This is also verified numerically in Table 1 and Table 2. 
\end{enumerate}
\begin{remark}\label{Remark.OrderAccuracyHypothesis}
Note that the restriction on $k$ ({\it i.e.} $k>1$) in the hypotheses of the above theorem is only to make $\delta$ defined by \eqref{Delta.exp} to tend to zero as $n_q\rightarrow \infty$.  Otherwise, the proof of the above theorem can go on with any $k\ge 1$ as in Theorem \ref{thm:errestimate} with appropriate choices of $\bfm{n}_p$, $\bfm{n}_q$, and $n$. Therefore, Theorem \ref{errorestimate.thm} may be used for $k\ge 1$, if $\bfm{n}_p$ and $\bfm{n}_q$ are appropriately chosen and fixed.
\end{remark}
\section{Numerical Comparison}\label{numerical.comparison.sec}
There are two tasks to be addressed numerically.  One is to examine the performance of the PiPCT approximation proposed in the previous section, and the other one is to check the agreement of numerical order of accuracy with the theoretical result obtained in Theorem \ref{errorestimate.thm}.

For the first task, we give numerical evidence that the proposed PiPCT method (in subsection \ref{piecewise.implementation.ssec}) captures singularities of a function without a visible Gibbs phenomenon.   
\begin{example}{\rm 
Consider the piecewise smooth function
	\begin{equation}\label{discondeg.eq}
	f(x) = 
	\begin{cases}
	x^3, &\quad\text{if } x \in [-1,-0.4),\\
	x^2+1, &\quad\text{if } x \in [-0.4,0.4),\\
	1.16 - (x-0.4)^{\frac{1}{2}}, &\quad\text{if } x \in [0.4,1].
	\end{cases}
	\end{equation}
The function $f$ involves a jump discontinuity at $x=-0.4,$ whereas at $x=0.4$ the function is continuous but not differentiable (referred in this article as {\it point singularity}) as shown in Figure \ref{fig:PiPChe_discondeg}(a). We fix the number of quadrature points as $n=200$ to reduce the error due to the quadrature formula considerably.  

The aim is to study the accuracy of the PiPCT approximant in terms of its parameters.  The aim is also to compare the performance of the PiPCT approximant with recently proposed global Pad\'e-based algorithms for approximating functions with singularities, namely, the singular Pad\'e-Chebyshev (SPC) method (see Driscoll and Fornberg \cite{dri-for_01a}, and Tampos {\it et al.} \cite{tam-etal_12a}), the robust Pad\'e-Chebyshev type (RPCT) method (see Gonnet \cite{gon-etal_11a, gon-etal_13a}), and finally with the well-known (global) PCT method.

For PiPCT approximant of $f$, we have chosen $N = 512$, $n = 200$, and $n_p^j=n_q^j=20$, for $j=0,1,\ldots,N-1,$ in all the uniformly discretized subintervals of $[-1,1]$. For all global algorithms, we have fixed the parameters as $n_p = 20 = n_q$, $n = 512 \times 200$, and $N = 1$.  
We can see from Figure \ref{fig:PiPChe_discondeg}(a) that in the smooth region, these approximants are well in agreement with the exact function, but in the vicinity of the singularities $x=-0.4$ and $x=0.4$, SPC could capture the jump discontinuity but not the point singularity (see the zoomed boxes).  Also, it is evident from the figure that the RPCT and the global PCT approximants are almost the same in a small neighborhood of the singularities.  From this figure, we can observe a significant role of piecewise implementation of the PCT approximation in capturing the singularities without using the locations of the singularities, unlike in the case of SPC.\vspace{0.1in}

\noindent {\it Numerical discussions as $N$ varies:} Figure \ref{fig:PiPChe_discondeg}(b) depicts the peaks of the pointwise error (as explained in Driscoll and Fornberg \cite{dri-for_01a}) for $N=2^k,$ $k=1,3,5,7,8,$ and $9.$ Figure \ref{fig:PiPChe_discondeg}(c) depicts the maximum error in the vicinity of $x=-0.4$ (in $-\!\!\!-\!\!\!$o$\!\!\!-\!\!\!-$ symbol) and in the vicinity of $x=0.4$ (in $-\!\!\!-\!\!\!*\!\!\!-\!\!\!-$ symbol) with logarithmic scale in the $y$-axis.
\begin{figure}[t]
	\centering
\includegraphics[height=14.cm,width=12cm]{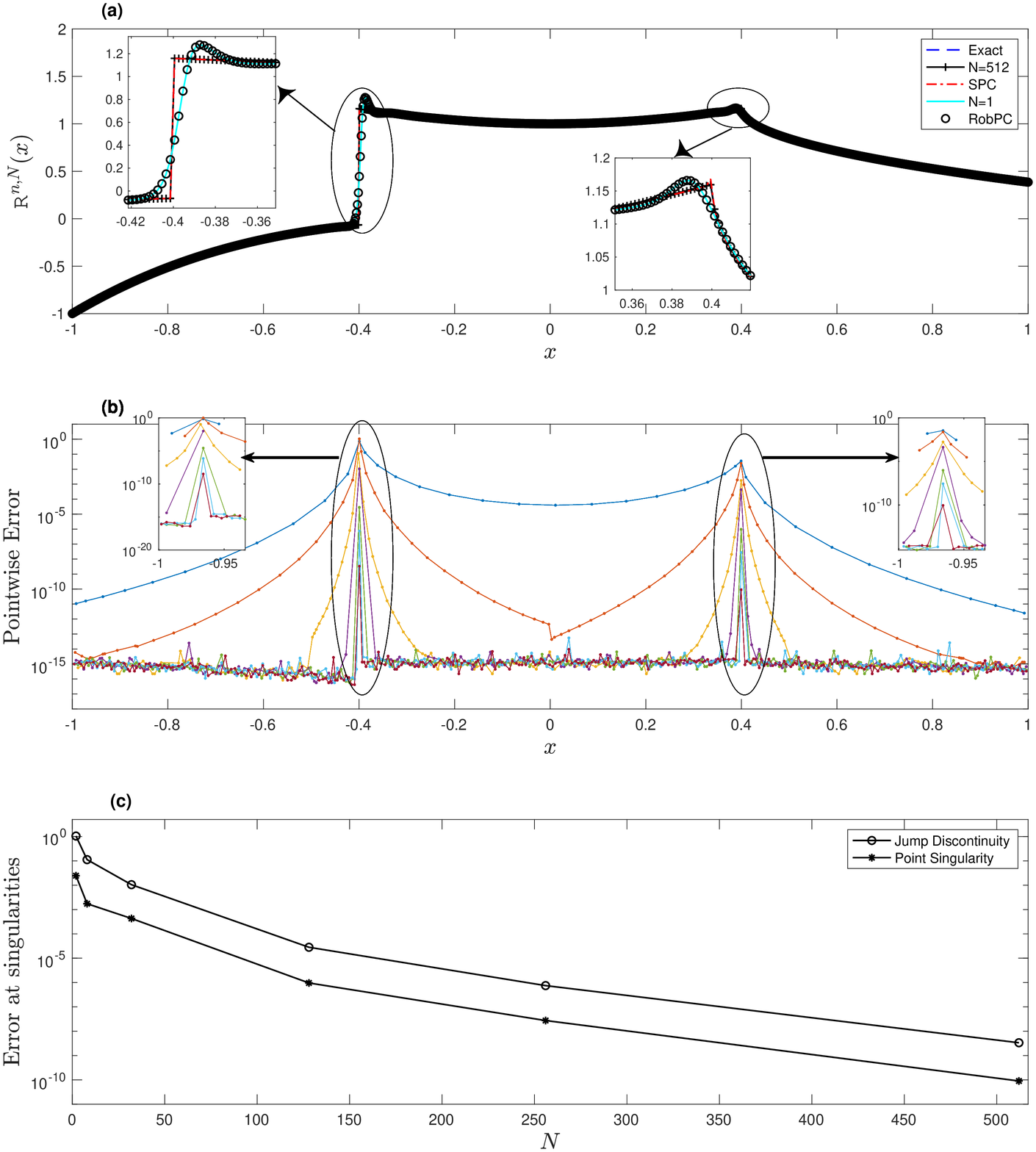}\vspace{-0.0in}
\caption{(a) Depicts comparison of PiPCT approximant of $f$ given by  \eqref{discondeg.eq} with the global Pad\'e-Chebyshev based algorithms, (b) depicts the peaks of pointwise error in approximating $f$ by PiPCT for $N = 2^k, k = 1,3,5,7,8,9$, and (c) depicts the $L^\infty$-error  in approximating $f$ by PiPCT  in the vicinity of $x = -0.4$ and $x = 0.4$.}
\label{fig:PiPChe_discondeg}
\end{figure}

Note that, $\mathsf{R}^{102400,1}(x)$ is a global approximation and $\mathsf{R}^{200,512}(x)$ is a piecewise approximation.   In both cases, we have given the values of the function at 102400 points.  In the case of the global PCT approximant, these points are the Chebyshev points in the interval $[-1,1]$, and we can observe an oscillation near the discontinuity, as shown in the zoomed box of Figure \ref{fig:PiPChe_discondeg}(a) in the vicinity of $x=-0.4$.  We also observe the same kind of behavior in the vicinity of the point singularity at $x=0.4$.  Note that the error estimate obtained in Theorem \ref{errorestimate.thm} shows that the sequence of piecewise approximants converges for functions that are at least continuous.  Though we do not have a theorem that gives convergence for functions involving jump discontinuities, the piecewise approximant in this example captures the singularities accurately, including the jump discontinuity at $x=-0.4$. Moreover, it is evident from Figure \ref{fig:PiPChe_discondeg}(b) and (c) that the sequence of piecewise approximants (in this example) tends to converge to the exact function as $N\rightarrow \infty$. 

Table \ref{ErrorNVaries.tab} compares the $L^1$-error and the numerical order of accuracy of the piecewise Chebyshev and the PiPCT approximants in the interval $[0.2,1]$ where the function $f$ has a point singularity at $x=0.4$.  Here, we can see the obvious advantage of using Pad\'e-Chebyshev type approximation when compared to the Chebyshev approximation. Also, we observe that the numerical order of accuracy is well in agreement with the theoretical result given in Theorem \ref{errorestimate.thm} with $k=1,$ $d^j=40$, for $j=0,1,\ldots,N-1$, and $n=200$ (see Remark \ref{Remark.OrderAccuracyHypothesis}).
\begin{table}
\caption{The $L^1$-error and the numerical order of accuracy of the piecewise Chebyshev and the PiPCT approximants of the function \eqref{discondeg.eq} in the interval $[0.2,1]$ as $N$ varies. We have taken $n_p=n_q=20$ in all the subintervals and $n=200.$}\vspace{0.1in}
\centering
\begin{tabular}{|l| c|c| c| c| c|  r|}
\hline
\multicolumn{1}{|c|}{$N$}&\multicolumn{2}{|c|}{Piecewise Chebyshev}&\multicolumn{2}{|c|}{Piecewise Pad\'e-Chebyshev type}\\
\hline
&$L^1$-error & order & $L^1$-error & order\\ 
\hline
2  &  0.603860 & --     &0.032616 &--     \\
\hline
8  &  0.03022378662560210039 & 2.725905     &0.00064588620006190815 &3.569908     \\
\hline
32  &  0.00109697967811300638 & 2.552230     &0.00002635315776778789 &2.462154     \\
\hline
128  &  0.00003324170968713331 & 2.564731     &0.00000001505864286582 &5.477419     \\
\hline
256  &  0.00000323501556305392 & 3.380088     &0.00000000021392558412 &6.171918     \\
\hline
512  &  0.00000004032820231387 & 6.343662     &0.00000000000035272088 &9.270412     \\
\hline
\end{tabular}
\label{ErrorNVaries.tab}
\end{table}
\begin{figure}[t]
	\centering
\includegraphics[height=13.cm,width=12cm]{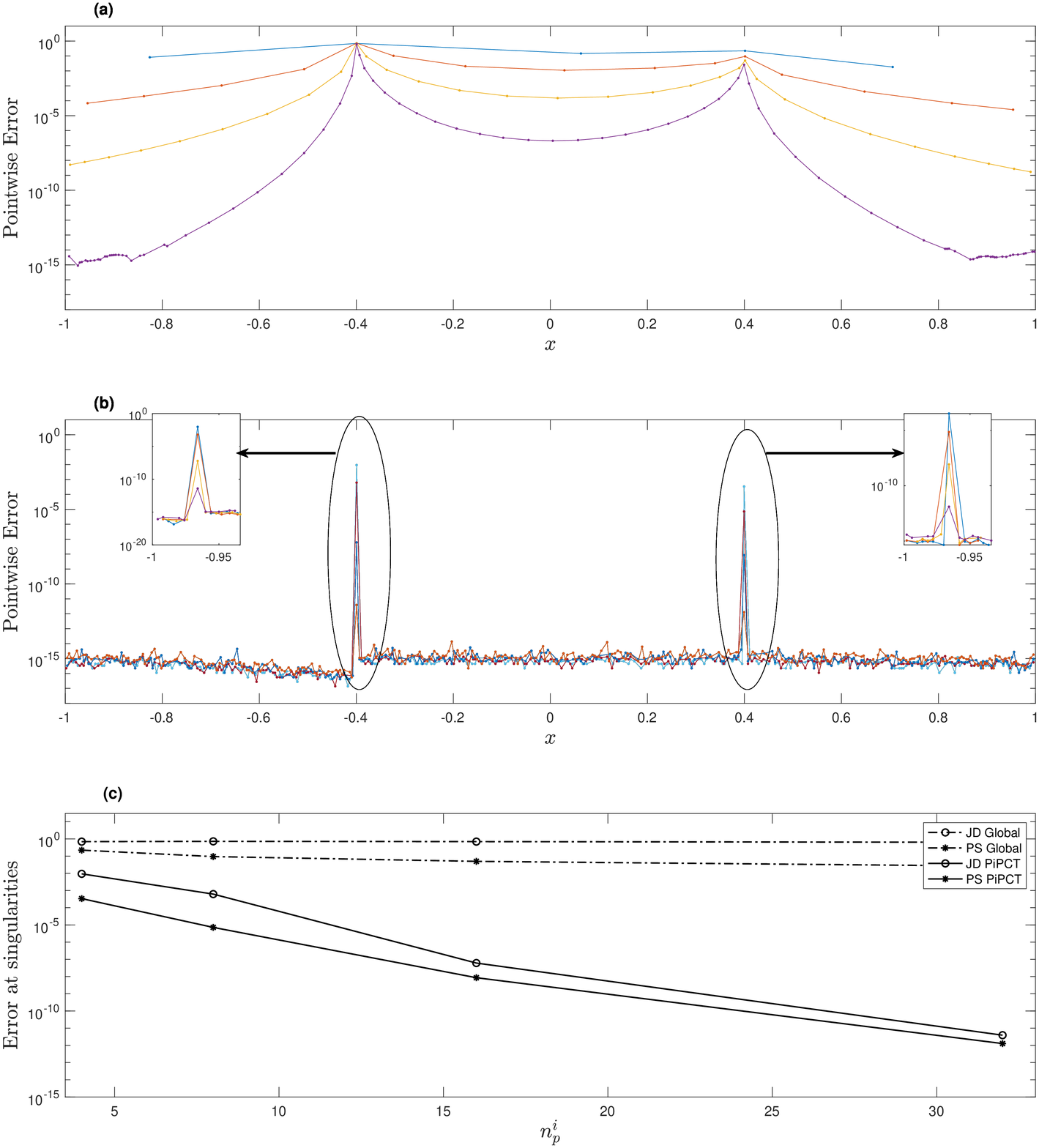}
\caption{Comparison between the global PCT and the PiPCT methods for $N = 512,$ $n = 200$ and $n_p = n_q = 2,8,16,32$. (a) depicts the peaks of the pointwise error $|f(x) - \mathsf{R}^{N\times n,1}_{\mbox{\tiny $\bfm{n}$}_p,\mbox{\tiny $\bfm{n}$}_q}(x)|$ (global PCT), (b) depicts the peaks of the pointwise error $|f(x) - \mathsf{R}^{n,N}_{\mbox{\tiny $\bfm{n}$}_p,\mbox{\tiny $\bfm{n}$}_q}(x)|$ (PiPCT), and (c) depicts the convergence of both the methods in the vicinity of $ x = -0.4$ (JD) and $x = 0.4$ (PS).}
\label{fig:PiPCheandChe_discondeg}
\end{figure}

\noindent {\it Numerical discussions as $n_p (=n_q)$ varies:}
Let us take $n_p^j=n_q^j=n_p$, $j=0,1,\ldots,N-1,$ and study the numerical convergence of the sequence of PiPCT approximantions $\mathsf{R}^{n,N}_{n_p,n_p}$ of the function $f,$ given by \eqref{discondeg.eq}, as $n_p$ varies. In order to examine the advantage of using the piecewise approximation, 
we also study numerically the convergence of the global Pad\'e-Chebyshev type approximantions $\mathsf{R}^{n\times N,1}_{n_p,n_p}$, for $n_p= 2,8,16,32$.

The peaks of the pointwise errors $|f(x)-\mathsf{R}^{n\times N,1}_{n_p,n_p}(x)|$ and $|f(x)-\mathsf{R}^{n,N}_{n_p,n_p}(x)|$ are depicted in Figures \ref{fig:PiPCheandChe_discondeg} (a) and (b), respectively. Here we have taken $n=200$, $N=512$ (therefore, in both the cases, we use the function values at 102400 grid points) and varied $n_p = 2, 8, 16$, and $32.$ We clearly observe that the pointwise error in the global implementation of the Pad\'e-Chebyshev type approximation does not seems to be converging at the jump discontinuity (at $x=-0.4$) as $n_p$ increases and the convergence is slow at the point singularity at $x=0.4$. This behaviour is more apparent in Figure \ref{fig:PiPCheandChe_discondeg} (c). On the other hand, from Figure \ref{fig:PiPCheandChe_discondeg} (b) and (c), we see that the convergence in the piecewise approximation is significantly faster both at the jump discontinuity and at the point singularity as $n_p$ increases.
\qed
}
\end{example}
In the above example, we have seen that the numerical order of accuracy is well in agreement with the theoretical result given in Theorem \ref{errorestimate.thm} when $k=1$ ({\it i.e.} when $k$ is odd).  
 
 In the following example, we show that this result also holds numerically for $k=2$.
 \begin{example}\label{SmoothPPC.ex}{\rm
 Consider the $C^1$-function
 $f(x) = x|x|$ 
 for $x\in [-1,1]$.  We use both the piecewise Chebyshev and the PiPCT approximations. The $L^1$-error and the order of accuracy are tabulated in Table \ref{SmoothPPC.tab}. The results are well in agreement with the theoretical results of Theorem \ref{errorestimate.thm} for $k=2$.
 }
 \qed 
 \end{example}
\begin{table}
\caption{The $L^1$-error and the numerical order of accuracy of the piecewise Chebyshev and the PiPCT approximants for the function given in Example \ref{SmoothPPC.ex}, as $N$ varies. We have taken $n_p=n_q=20$ in all the subintervals and $n=200.$}\vspace{0.1in}
\centering
\begin{tabular}{|l| c|c| c| c| c|  r|}
\hline
\multicolumn{1}{|c|}{$N$}&\multicolumn{2}{|c|}{Piecewise Chebyshev}&\multicolumn{2}{|c|}{Piecewise Pad\'e-Chebyshev}\\
\hline
&$L^1$-error & order & $L^1$-error & order\\ 
\hline
2  &  0.00000000000001916544 & --     &0.00000000000002741904 &--     \\
\hline
4  &  0.00000000000000282010 & 3.751452     &0.00000000000000335724 &4.111221     \\
\hline
8  &  0.00000000000000028910 & 3.875159     &0.00000000000000031289 &4.037213     \\
\hline
16  &  0.00000000000000003397 & 3.366981     &0.00000000000000003508 &3.440576     \\
\hline
\end{tabular}
\label{SmoothPPC.tab}
\end{table}
\section{Construction of an Adaptive Partition}\label{APiPCT.sec}
In section \ref{numerical.comparison.sec}, we demonstrated numerically the performance of the PiPCT approximation in capturing singularities of a function accurately.
 However, the numerical results depicted in Figure \ref{fig:PiPChe_discondeg} suggest that we need a sufficiently finer discretization to get good accuracy in the vicinity of the singularities.  Such a finer discretization is not needed in the regions where the function is sufficiently smooth.  
 This motivates us to look for an adaptive implementation of the PiPCT method.

The main idea of our adaptive algorithm is to identify the subintervals (referred to as {\it bad-cells}) where the function has a singularity and then bisect the subinterval.  Therefore, as a first step, we need a strategy to identify the bad-cells.
\subsection{Singularity Indicator}
In general, rational (Pad\'e) approximations are better than (see for instance, \cite{kab-mad_05a,bak-pet_82a}) the polynomial approximations in the case of approximating non-smooth functions. A natural concern about Pad\'e approximations is the poles of the approximants. An approximant may have poles at places where the function has no singularities. In literature, such poles are called {\it spurious poles}  (for precise definition, see \cite{sta_98a}). Besides such spurious poles, the poles of a PCT approximant of a non-smooth function also lie in a sufficiently small neighborhood of singularities. Poles are intractable on a computer because of the presence of a round-off error. Nevertheless, the following result can be observed numerically:
\begin{result}\label{res:jumplocate}
Let $f\in L^1[-1,1]$, which is at most continuous, with an isolated singularity at $x_0 \in (-1,1)$ and let $n>0$ be a given integer. For a given  $\epsilon>0$, there exists an integer $n_{p_0}>0$,  such that
\begin{equation}\label{Adaptive.Condition.eq}
|Q_{n_{p}}^n(z)|<\epsilon,~~\mbox{\rm for all }~~n_p\ge n_{p_0},
\end{equation}
whenever $|x_0-{\rm Re}(z)|<\delta,$ for some $\delta=\delta(n_p,\epsilon)>0.$ 
Here, $Q_{n_p}^n(z)$ $($for all $z\in \mathbb{C}$ such that $|z|=1)$ is the denominator polynomial of the PCT approximant $\mathsf{R}_{n_p,n_p}^n$ of $f$.
\end{result}
We illustrate this result in the following numerical example.
\begin{figure}[t]
	\centering
\includegraphics[height=5.cm,width=12cm]{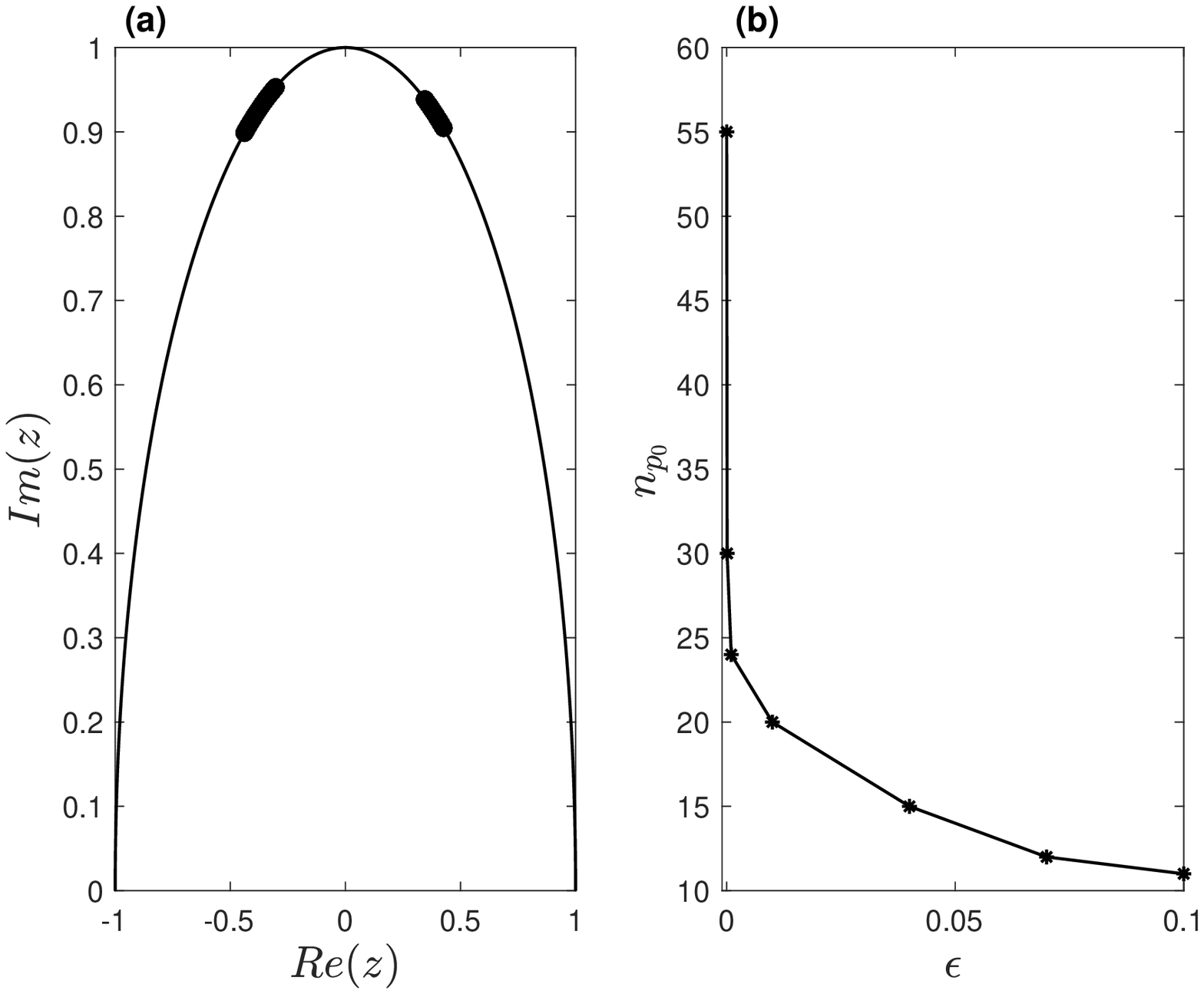}
\caption{(a) Depicts the points (with `o' symbol) on the unit circle which satisfies $|Q_{20}(z)| < 10^{-2}$, and (b) depicts the $\epsilon$ values and the corresponding values of $n_{p_0}$ that satisfy \eqref{Adaptive.Condition.eq}.}
\label{fig:Singularity.Location.discondeg}
\end{figure}
\begin{example}
{\rm 
Consider the piecewise smooth function $f$ given by \eqref{discondeg.eq}, where the function has a jump discontinuity at $x=-0.4$ and a point singularity at $x=0.4.$ Let $n=200.$ Figure \ref{fig:Singularity.Location.discondeg} (a) depicts the points (`o' symbol) at which $|Q_{20}^{200}(z)|\le 10^{-2}.$ We observe that the real parts of these points are accumulated near  $x=-0.4$ and $x=0.4$.   Further, Figure \ref{fig:Singularity.Location.discondeg} (b) depicts the graph whose $x$-coordinate represents different values of $\epsilon$ and the $y$-coordinate is taken to be the corresponding values of $n_{p_0}$, which are the minimum degrees of the denominator polynomial $Q_{n_p}^{200}$, for which the condition \eqref{Adaptive.Condition.eq} holds numerically. 
\qed
}
\end{example}
Based on the above numerical observation, we define a notion of  {\it bad-cells} in a given partition as follows:
\begin{definition}\label{eps.bad.cell.def}{\rm
Let $f\in L^1[a,b]$ and let $P_N$ be a given partition. For a given $\epsilon>0$ and the positive integers $m$ and $n$, a subinterval $I_j=[a_j,b_j]$ with $a_j\in P_N$ is said to be an {\it $\epsilon$-badcell} if  $|Q_{m}^n(z)|<\epsilon$ for some $z\in \mathbb{C}$ on the unit circle with ${\rm Re}(z)\in \mathsf{G}_j^{-1}(I_j)$.
}
\end{definition}
\subsection{Generation of an Adaptive Partition}\label{generation.adaptive.partition.ssec}
In this subsection, we propose an algorithm to generate a partition that consists of finer discretization in the vicinity of singularities.

Consider a function $f\in L^1[a,b]$.  Choose an $\epsilon>0$ sufficiently small, a tolerance parameter $\tau>0$, sufficiently large positive integer $n$, and $m$ $(\ll n/2)$ sufficiently small with
 $\bfm{n}_p=\bfm{n}_q=(m,m).$ 

Let us take $N_0=2$, $a_0=a$ and $b_0=b.$ Let $I^{(0)}_0$ and $I^{(0)}_1$ be the subintervals of the partition $P_{N_0}:=\{a_0,(a_0+b_0)/2,b_0\}.$ 

Perform the following steps for $k=0,1,2,\ldots$:
\begin{enumerate}
\item For each $j=0,1,\ldots, N_k-1$, obtain the denominator polynomial, denoted by $Q_{m,j}^n(z),$  of the Pad\'e-Chebyshev type approximant of the function $f|_{I^{(k)}_j}$.
\item For each $j$, check if the subinterval $I^{(k)}_j$ is an $\epsilon$-badcell as per the condition given in Definition \ref{eps.bad.cell.def}.
\item If a subinterval $I^{(k)}_j$ is detected as an $\epsilon$-badcell, then bisect this subinterval and add the end points to the partition $P^{b}_{N_{k+1}}$, referred to as the {\it badcells partition}.
\item Define the new partition $P_{N_{k+1}} = P_{N_{k}} \cup P^{b}_{N_{k+1}}$ and denote the subintervals of this partition as $I^{(k+1)}_j$, for $j=0,1,2,\cdots, N_{k+1}-1.$
\item Let $l_{*}: = \min\{|I^{(k+1)}_j|~|~j=0,1,2,\cdots, N_{k+1}-1\}$. If $l_*<\tau,$ then stop the iteration. Otherwise, replace $k$ by $k+1$ and repeat the above process (note that it is enough to repeat the process only for those intervals $I^{(k)}_j$ which are identified as $\epsilon$-badcells).
\end{enumerate}
The final outcome of the above process is the required {\it adaptive partition}, which we denote by $P^{*}_{N_K}$.
\subsection{Degree Adaptation}\label{degree.adaptation.subs}
We expect the adaptive partition proposed above to give an efficient way of obtaining the PiPCT approximation of a function $f\in L^{1}[a,b]$ with isolated singularities. In the generation of the adaptive partition, we fixed $n_p=n_q$ in all the subintervals of the partition.  The numerical results depicted in Figure \ref{fig:PiPCheandChe_discondeg}(c) suggest that we can improve the accuracy by choosing higher degrees for the numerator and denominator polynomials in the $\epsilon$-badcells of the partition  $P^{*}_{N_K}$.  However, as mentioned in Remark \ref{UB.behavior.remark},
the upper bound decreases as $d$ increases with $d<n$, whereas the opposite behavior of the upper bound is apparent when $d>n$. This suggests us to choose $m< n/2$ sufficiently small in the smooth regions.  

In this subsection, we propose a better choice of the degree of the numerator polynomial in the $\epsilon$-badcells. 

\begin{figure}
\centering\hspace{-0.2in}
\includegraphics[height=5cm,width=12cm]{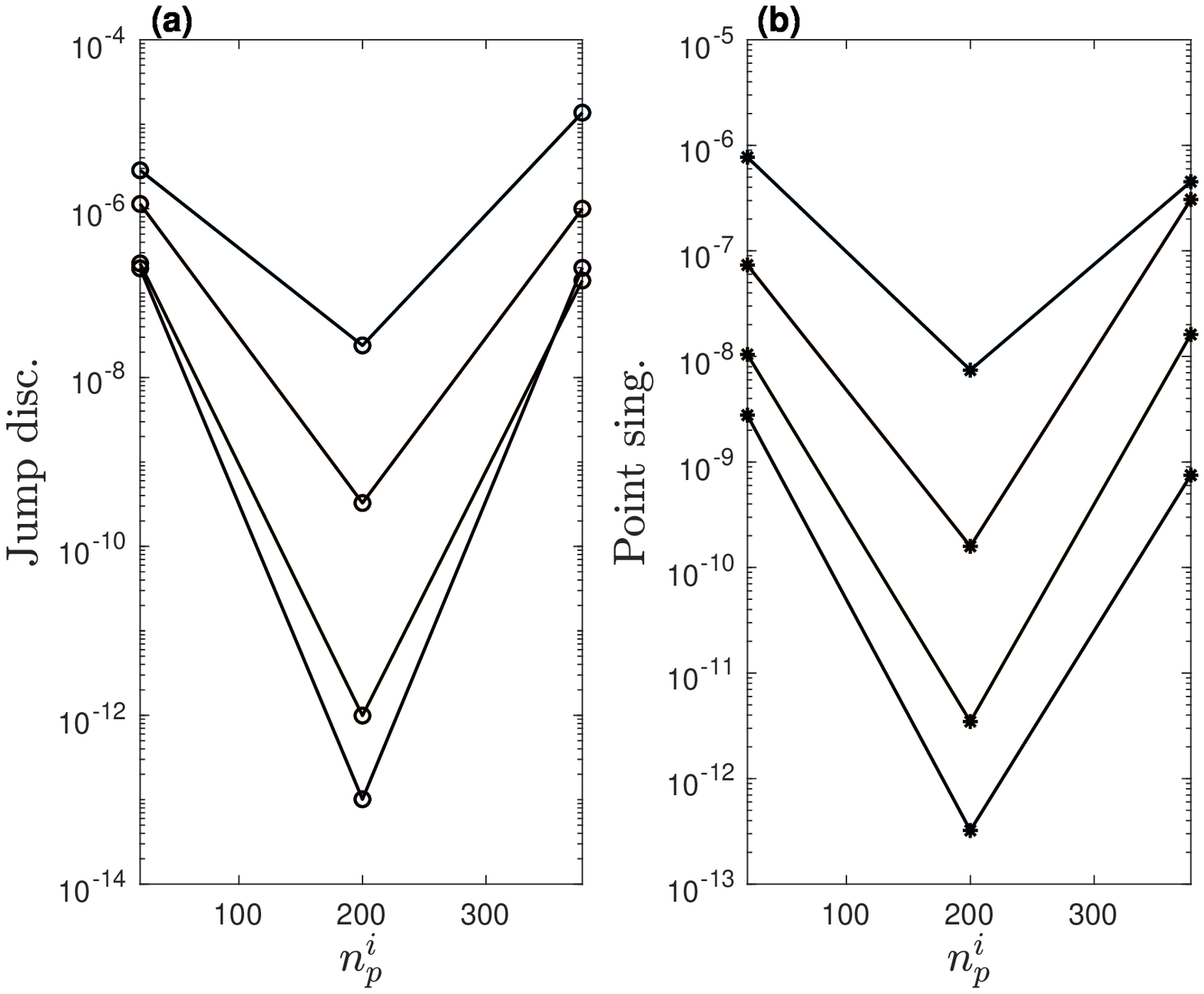}	
\caption{Pointwise error in PiPCT approximation of the function $f$ given by  \eqref{discondeg.eq} in a small neighborhood of \textbf{(a)} the jump singularity at $x = -0.4$ and \textbf{(b)} the point singularity at $x = 0.4.$ Different lines (from top to bottom) correspond to $N = 104, 208, 312, 416$, where $n_q^i=m=20$ and $n=200$. Three points lying on a graph (shown in symbol) correspond to the case when (from left to right) $n_p = n_q,$ $n_p=n,$ and  $n_p=2n-n_q-1$, respectively.}
\label{fig:Adaptivetool}
\end{figure}

From Proposition \ref{prop:1}, we see that the coefficient vector $\bfm{q}$ of the denominator polynomial of $\mathsf{R}_{m,m}^n$ is the flip vector of the denominator polynomial of $\mathsf{R}_{2n-m-1,m}^n$  for a given $m<n-1.$  Since the coefficient vector $\bfm{p}$ of the numerator polynomial of a Pad\'e-Chebyshev type approximant is computed using $\bfm{q}$, we expect that the errors involved in approximating $f$ in an $\epsilon$-badcell by $\mathsf{R}_{m,m}^n$ and $\mathsf{R}_{2n-m-1,m}^n$ are almost the same. 
In particular, this result has been verified numerically for the function $f$ given by \eqref{discondeg.eq} when $m$ is sufficiently smaller than $n-1$ and is shown in Figure \ref{fig:Adaptivetool}. This figure suggests that the error attains its minimum when $n_p=n$ in this example.

Based on the above numerical observation, we propose the following choice of the polynomial degrees in the adaptive partition obtained in subsection \ref{generation.adaptive.partition.ssec}. Note that this is only a preliminary idea for degree adaptation.  A rigorous mathematical analysis is needed to obtain a more reliable procedure.
\subsubsection{ Choice of polynomial degrees in the adaptive algorithm:}\label{subsubsec:choiceofdegree}
For a given number of quadrature points $n$, choose $m\ll n/2$.  Set $n_p=n$ and $n_q=m$ in all $\epsilon$-badcells of the adaptive partition $P^*_{N_K}$ and set $n_p=n_q=m$ for other cells of the partition. 

We call the resulting algorithm as {\it adaptive piecewise Pad\'e-Chebyshev type} algorithm (APiPCT).
\begin{figure}[t]
\centering
\includegraphics[height=12.5cm,width=12cm]{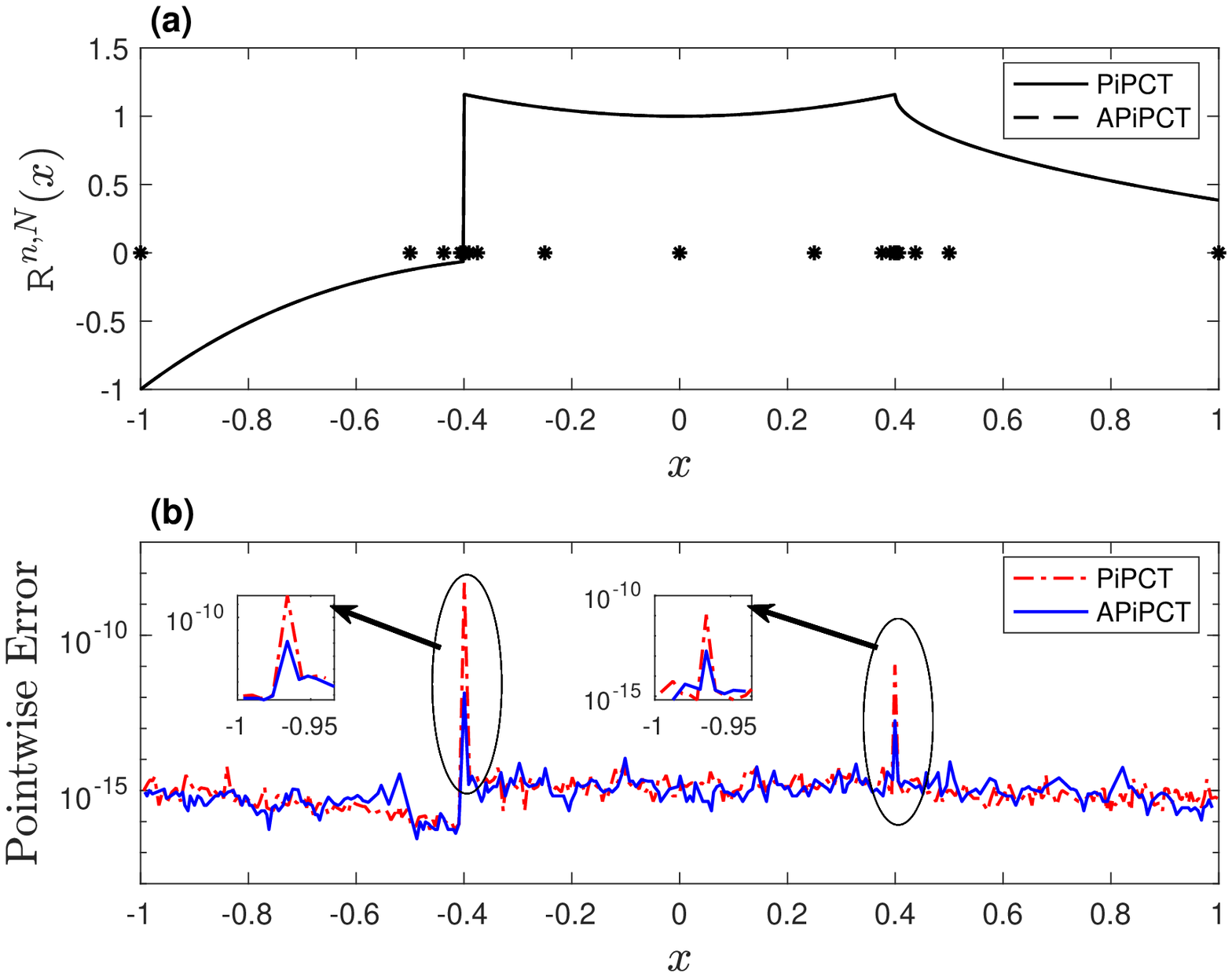}
\caption{\textbf{(a)} Depicts the comparison of the approximants obtained by the PiPCT algorithm and the APiPCT algorithm. The symbols `*' represent the points in the adaptive partition $P^*_{N_K}$. \textbf{(b)} Depicts the comparison of the pointwise error of the approximants using the two algorithms for $N = 512,n = 100$, and for numerator and denominator degrees $n_p = n_q = 20$ for the function $f$ given by \eqref{discondeg.eq}.}
\label{fig:plotErrPiPCAPiPC}
\end{figure}
\subsection{APiPCT Method}
We summarize the APiPCT algorithm as follows:

Consider a function $f \in L^1[a,b]$. We construct an APiPCT approximation of $f$ in the following steps:
\begin{enumerate}
\item Generate an adaptive partition $P^*_{N_K}$ as explained in section \ref{generation.adaptive.partition.ssec}.
\item \label{APiPCT2} Choose the degrees of numerator and denominator polynomials stated in subsection \ref{subsubsec:choiceofdegree}.
\item Use PiPCT algorithm given in section \ref{piecewise.implementation.ssec} with partition $P^*_{N_K}$ and degrees $\bfm{n}_p$ and $\bfm{n}_q$ as given in Step \ref{APiPCT2} above.
\end{enumerate}
The outcome of the above algorithm is a PiPCT approximation of $f$ with respect to the given partition $P^*_{N_K}$ and is given by
\begin{equation}
\mathsf{R}_{n_p,n_q}^{n,N_K}(t) =
\begin{cases}
\mathsf{R}_{n,m}^{n}(t) &\quad\text{if } t \in \mbox{an }\epsilon\mbox{-badcell}\\
\mathsf{R}_{m,m}^{n}(t) &\quad\text{otherwise}.
\end{cases}
\end{equation}
\section{Numerical Experiments Using Adaptive Partitioning}\label{numerical.experiment.APiPCT.sec}
To validate the APiPCT algorithm, we approximate the function $f$ given by \eqref{discondeg.eq} using PiPCT and  APiPCT algorithms, and compare the approximants in Figure \ref{fig:plotErrPiPCAPiPC}(a).  Here, we choose $n=100$, $m=20$ (common to both the algorithms), and $N=512$.  Further, for the APiPCT algorithm, we choose $\epsilon =  10^{-2}$ and $\tau=(b-a)/N=1/256.$  With these parameters, the adaptive algorithm generates 18 cells as depicted by `*' symbol in Figure \ref{fig:plotErrPiPCAPiPC}(a).  Figure \ref{fig:plotErrPiPCAPiPC} (b) depicts the corresponding pointwise errors.  Here, we observe that the error in the badcells is significantly reduced in the APiPCT algorithm when compared to that of PiPCT. This improvement is mainly because of the choice of numerator polynomial degree $n_p=n$ in the $\epsilon$-badcells.

\begin{figure}[t]
	\centering
	\includegraphics[height=5.5cm,width=13cm]{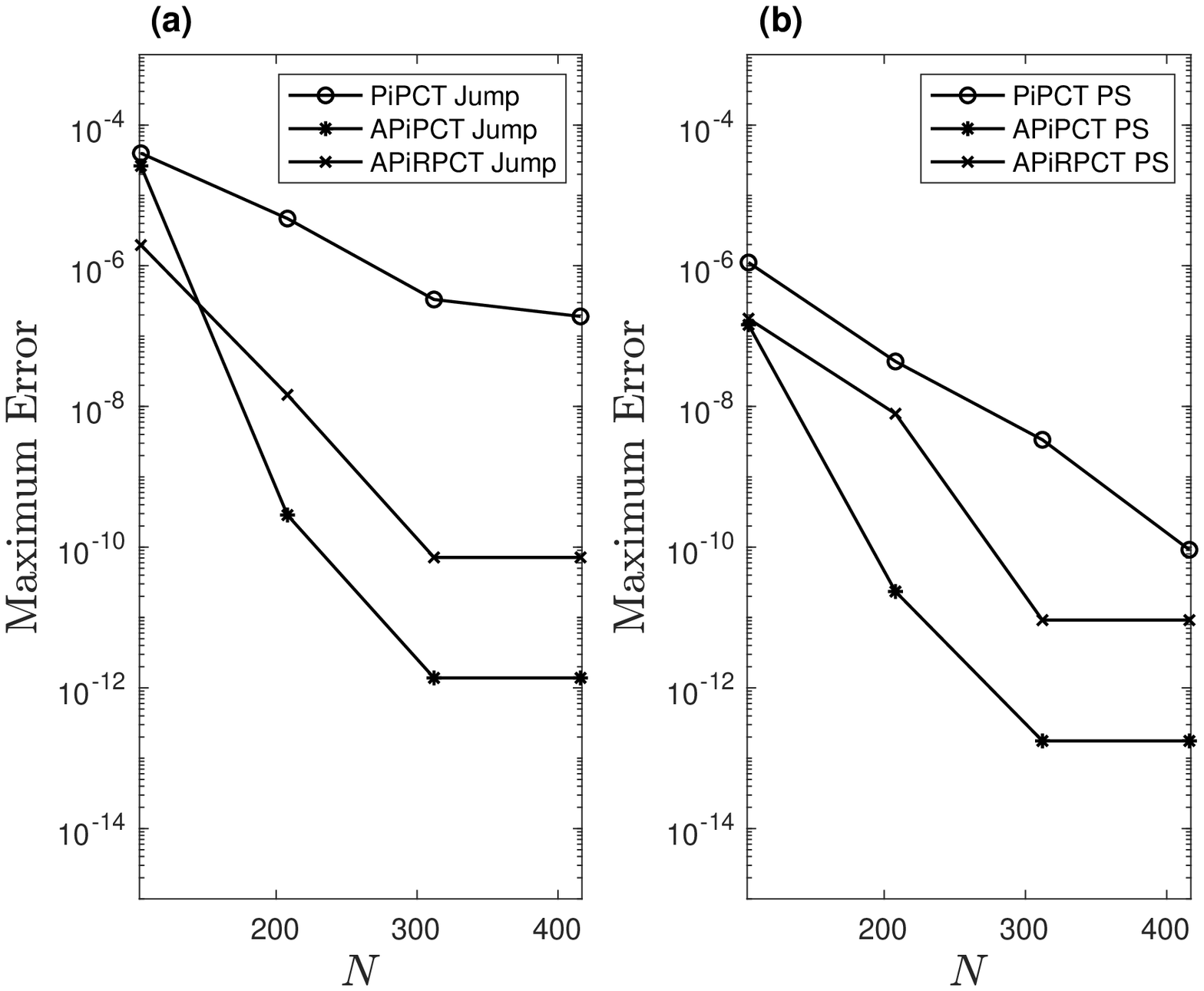}	
	\caption{Comparison between PiPCT, APiPCT and APiRPCT algorithms. \textbf{(a)} Depicts the comparison in maximum error in a neighborhood of the jump discontinuity at $x = -0.4$ and \textbf{(b)} depicts the comparison in maximum error in a neighborhood of point singularity at $x = 0.4$,  for $N = 104, 208, 312, 416$, and for numerator and denominator degrees $n_p = n_q = 20$ for function $f$ given by \eqref{discondeg.eq} on the domain $[-1,1]$. }
	\label{fig:Figure62}
\end{figure}

We are also interested in studying the performance of the robust Pad\'e-Chebyshev  (RPCT) method of Gonnet {\it et al.} \cite{gon-etal_13a} in the adaptive piecewise algorithm developed in the above section.  Observe that without any further modifications, we can replace the PCT method by the RPCT method  in the adaptive algorithm. We denote the adaptive piecewise RPCT algorithm by APiRPCT. We observe that the APiRPCT method also captures the singularities as sharply as the APiPCT algorithm (Figure \ref{fig:Figure62}).

Figures \ref{fig:Figure62} {(a)} and {(b)} depict the maximum error in the vicinity of singularities present in the function \eqref{discondeg.eq} at $x = -0.4$ (jump discontinuity) and $x = 0.4$ (point singularity), respectively, for $N =
104, 208, 312, 416$. We can see the performance of all the three methods (PiPCT, APiPCT, and APiRPCT) near a singularity and the significant role of choosing the numerator degree adaptively. Here, we observe that for a given $n_q$, setting $n_p=n$ in the adaptive methods decreases the maximum error in an $\epsilon$-badcell more rapidly than in PiPCT as $N$ increases. 

We also observe from Figures \ref{fig:Figure62} {(a)} and {(b)} that the maximum error in the APiPCT method decreases more rapidly than that of the APiRPCT method. This behavior is because the RPCT method calculates the rank of the Toeplitz matrix (using singular value decomposition) and reduces the numerator and denominator degrees diagonally to reach to the final position (correspond to the minimum degree denominator) in the upper left corner of a square block of the Pad\'e table (for block structure details see Gragg \cite{gra_72a} and Trefethen \cite{tre_84a}). Thus, in the process of minimizing the occurrence of spurious pole-zero pairs, the RPCT method decreases the denominator degree and therefore reduces the accuracy (as already been mentioned in Gonnet {\it et al.} \cite{gon-etal_13a}).  

\begin{figure}[t]
\centering
\includegraphics[height=5.5cm,width=8cm]{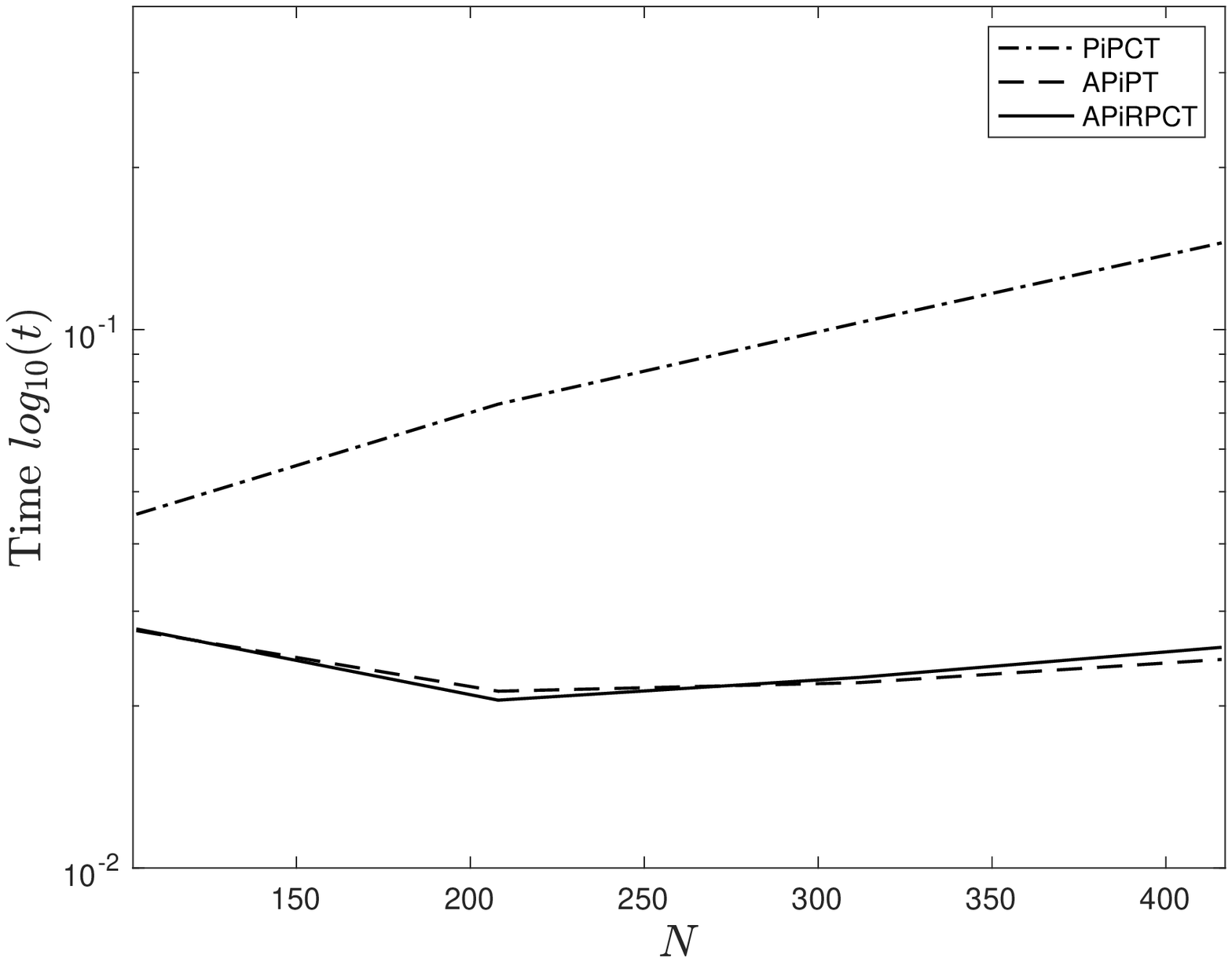}	
\caption{Comparison between the time taken by the PiPCT, the APiPCT and the APiRPCT algorithms to compute the approximants of the function $f$ given by \eqref{discondeg.eq}, where $N = 104, 208, 312, 416$, and the numerator and the denominator degrees are $n_p = n_q = 20$. }
\label{fig:Efficiency}
\end{figure}

Finally, to study the efficiency of the adaptive piecewise algorithm, we compare the time taken by PiPCT, APiPCT, and APiRPCT algorithms to approximate the function given by \eqref{discondeg.eq}. Figure \ref{fig:Efficiency} depicts the comparison between the time taken by these three algorithms as we increase the number of partitions $N$ (correspondingly decreasing the tolerance parameter $\tau=(b-a)/N$). In this figure, we observe that the time taken by PiPCT approximation increases with $N$ while the time taken by the adaptive algorithms APiPCT and APiRPCT remain almost the same. Although the RPCT approximation does a repeated singular value decomposition to get a full rank matrix, finally, the construction of the robust Pad\'e approximation is done possibly with a lower degree polynomials and hence the time taken in the repeated singular value decomposition is compensated in the Pad\'e construction.  Whereas, in the APiPCT algorithm,  the Pad\'e approximation is performed with higher-order polynomials. Hence the times taken by the APiPCT and the APiRPCT algorithms finally remain almost the same in this example.

\section{Comments on Froissart Doublets}\label{comment.FD.sec}
From a theoretical point of view, it is known that (see, for instance, Baker and Peter \cite{bak-pet_82a}) a Pad\'e approximation accelerates the convergence of a truncated series. It can be used as a noise filter in signal processing. It significantly reduces the effect of Gibbs oscillations but not able to eliminate it. Apart from these (theoretical) properties of a Pad\'e (rational) approximation, there are difficulties in elucidating the approximation power of these approximants correctly. It is mainly because of the random occurrence of the poles of the PCT approximant in the complex plane. For a real-valued function (after projecting onto the complex plane), if these poles are sufficiently away from the unit circle, then it may not affect the approximation. The problem occurs when an approximant has poles in a region where the function has no singularities. In such cases, one can not expect accurate result (near the poles). These poles are called \textit{spurious poles}. Baker \textit{et al.} \cite{bak-etal_61a} shows that spurious poles are isolated and always accompanied by a zero. The mentioned spurious pole-zero pair is also known as \textit{Froissart doublets}. The occurrence of Froissart doublets is a fundamental mathematical issue; these are the difficulties in establishing convergence theorems of Pad\'e approximants of order $[n_p/n_p]$, which cannot be achieved without a restriction of convergence in measure or capacity and not uniform convergence \cite{pom_73a,nut_70a}. On a computer in floating-point arithmetic they even arise more often and are known as \textit{numerical Froissart doublets} (see Nakatsukasa {\it et al.} \cite{nak-etal_18a}). 

There are at least two methods to recognize numerical Froissart doublets. In \cite{gon-etal_13a,gon-etal_11a,nak-etal_18a}, authors proposed a way to identify those doublets by calculating the absolute values of their residual. Another method to recognize the spurious pole-zero pair is by examining the distance between the two \cite{bes-per_09a}. Here, we use the residual method to recognize the spurious poles.  

\begin{figure}[]
	\centering
	\includegraphics[height=5.2cm,width=5cm]{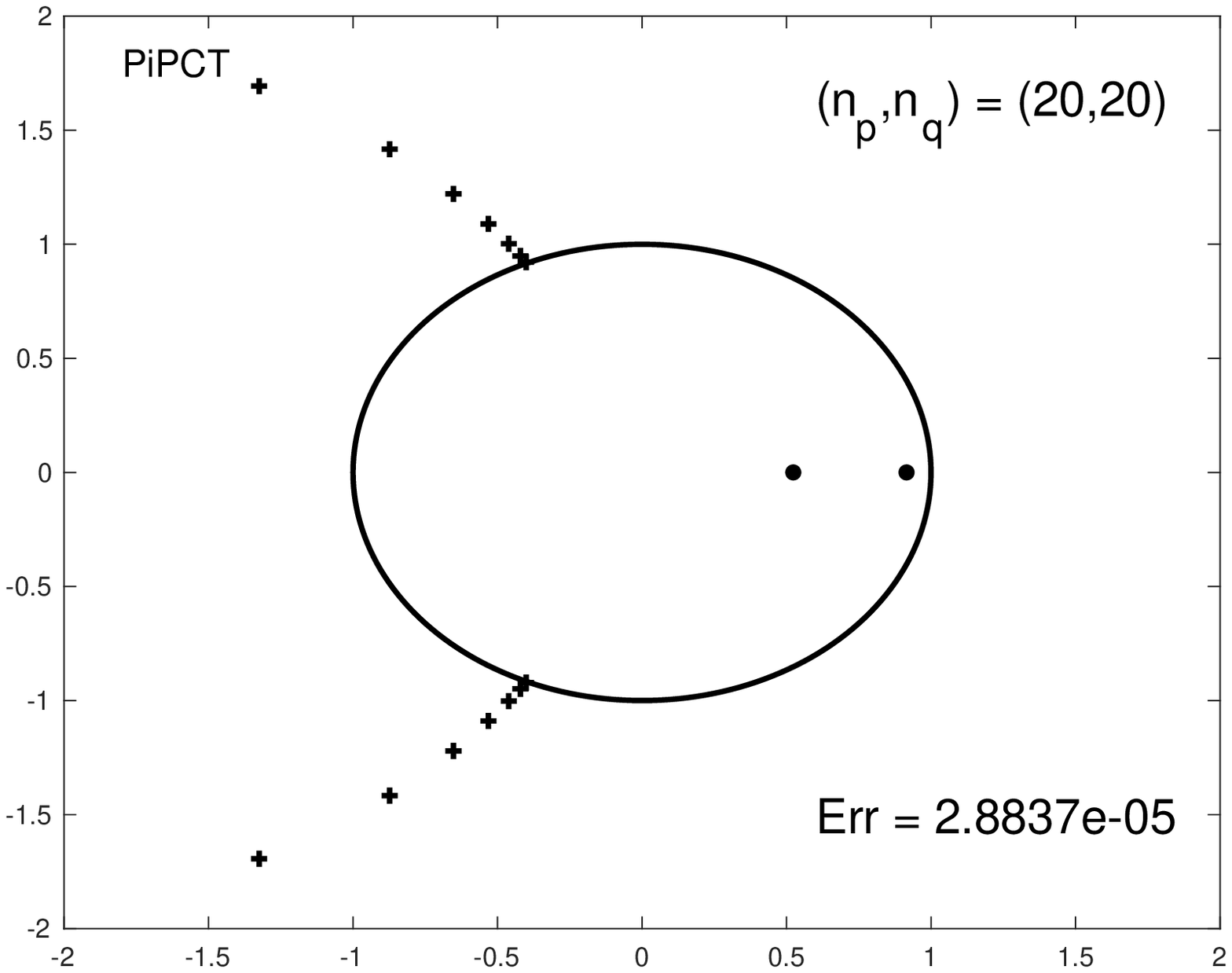}	
	\includegraphics[height=5.2cm,width=5cm]{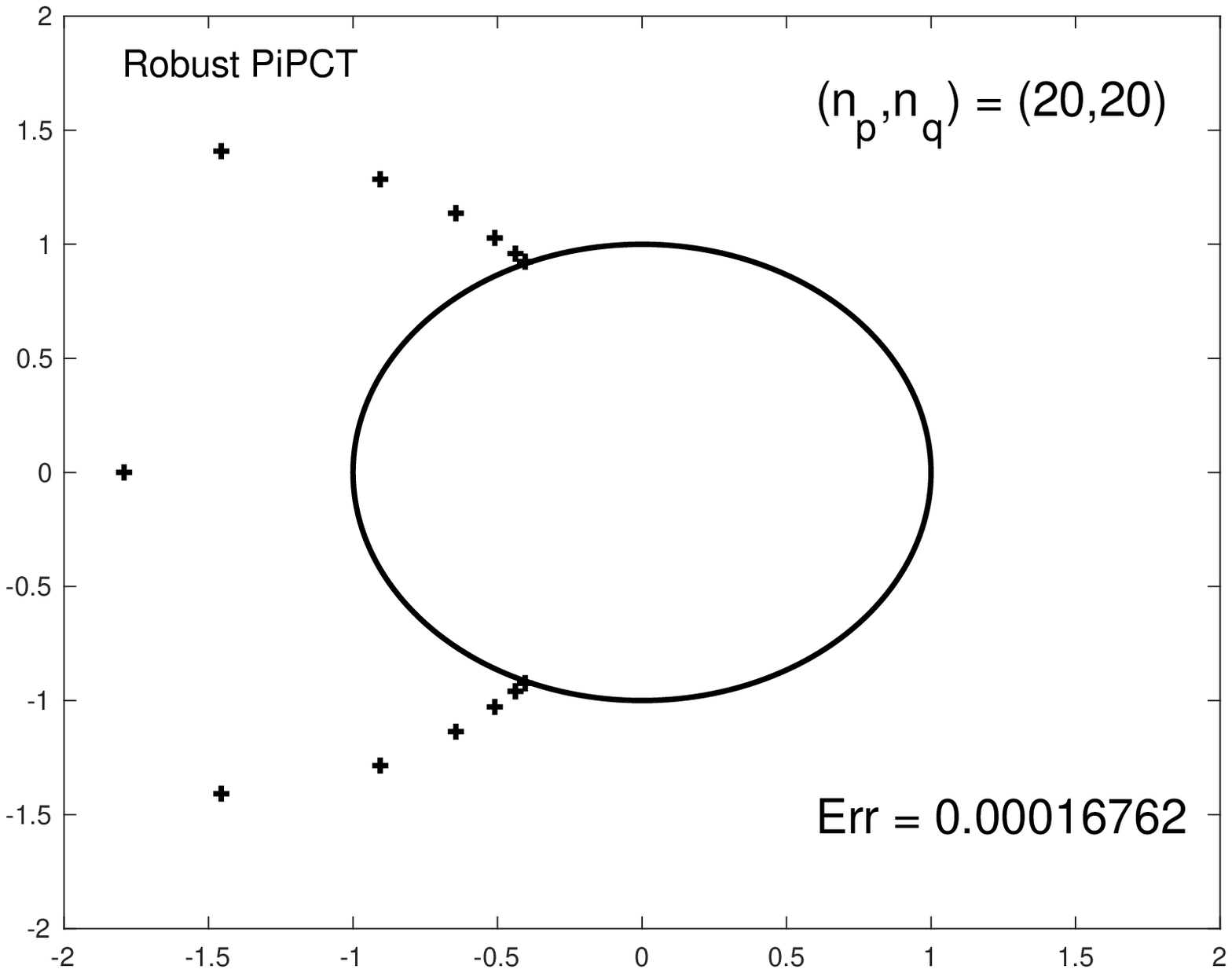}\\
	\includegraphics[height=5.2cm,width=5cm]{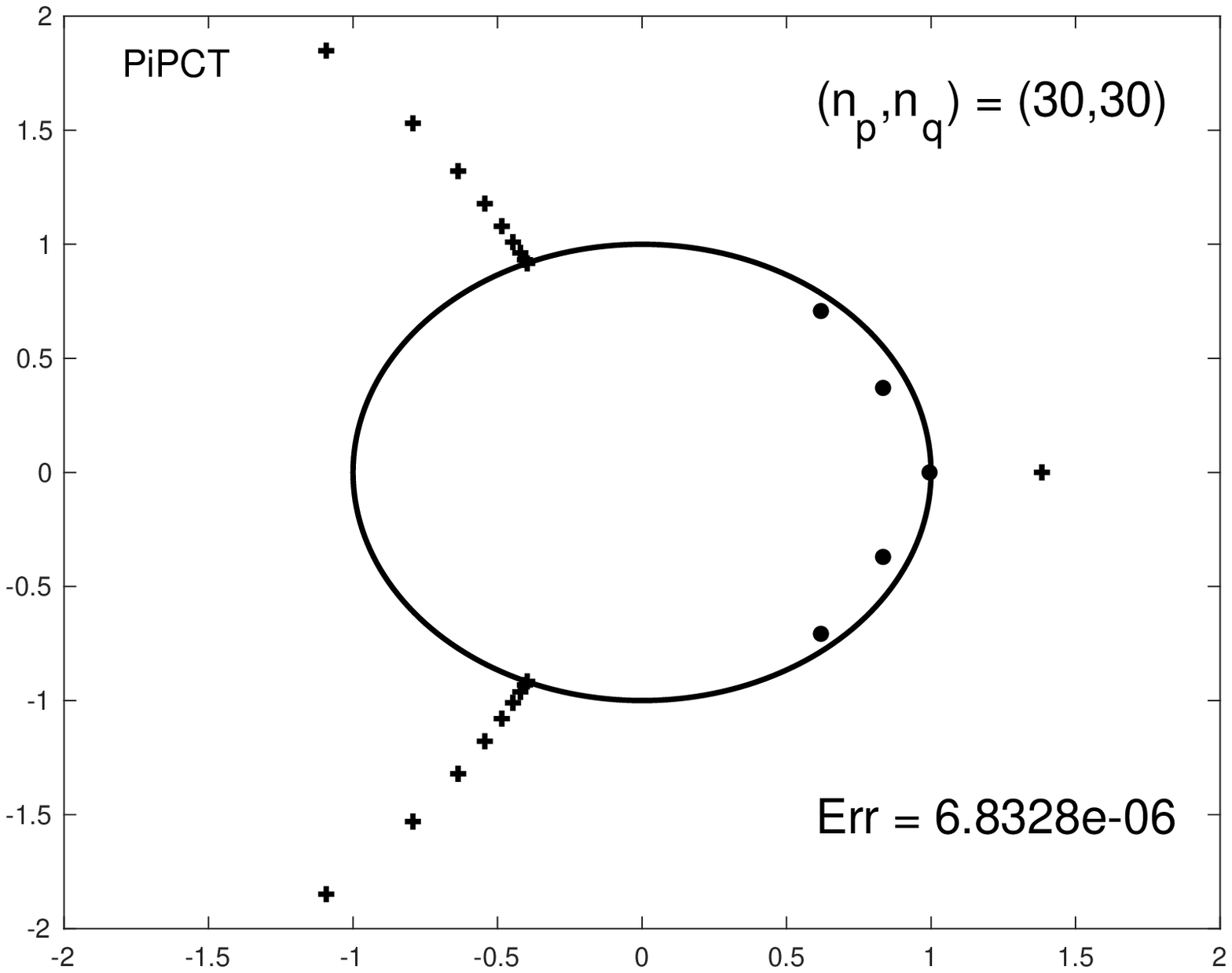}	
	\includegraphics[height=5.2cm,width=5cm]{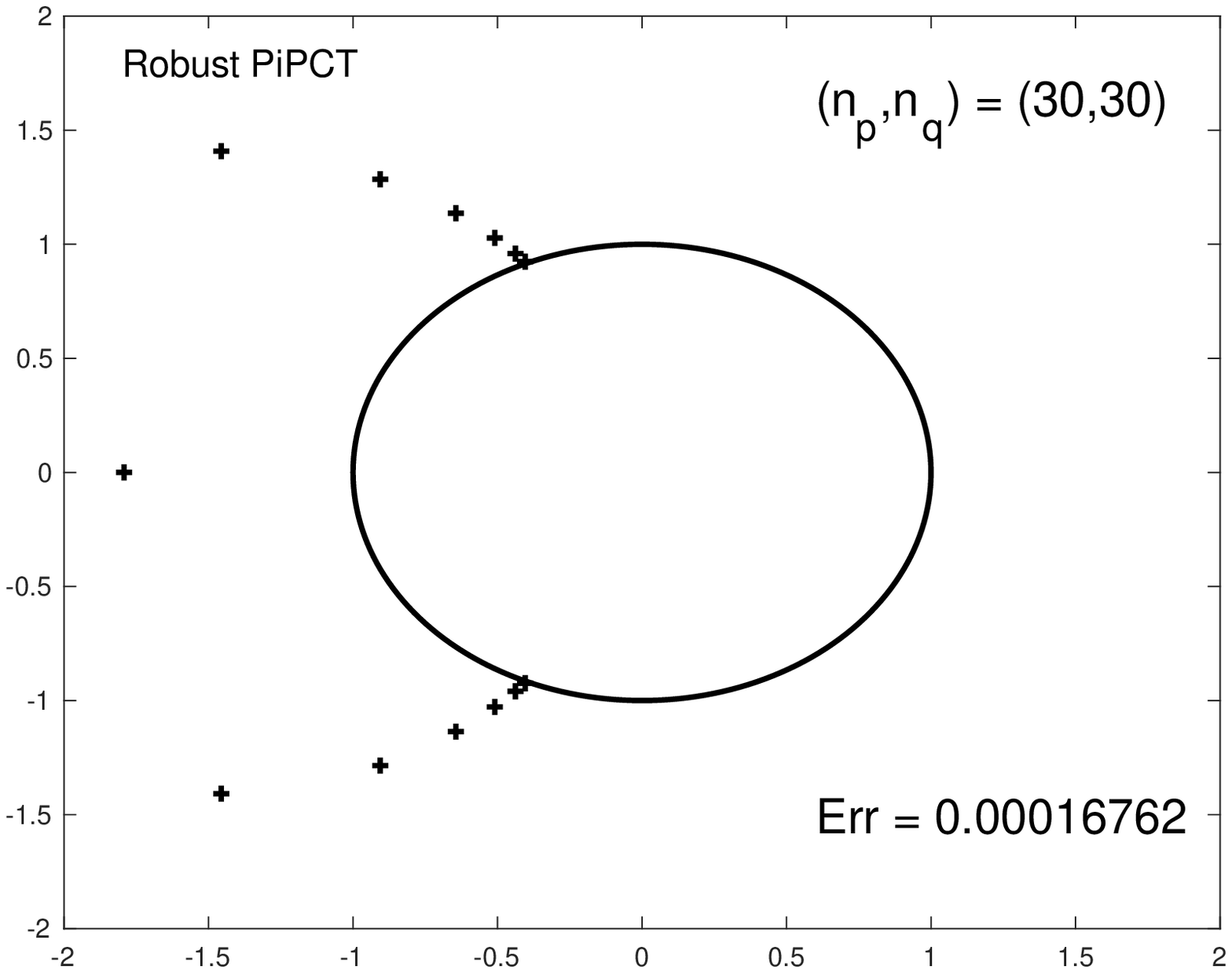}\\
	\includegraphics[height=5.2cm,width=5cm]{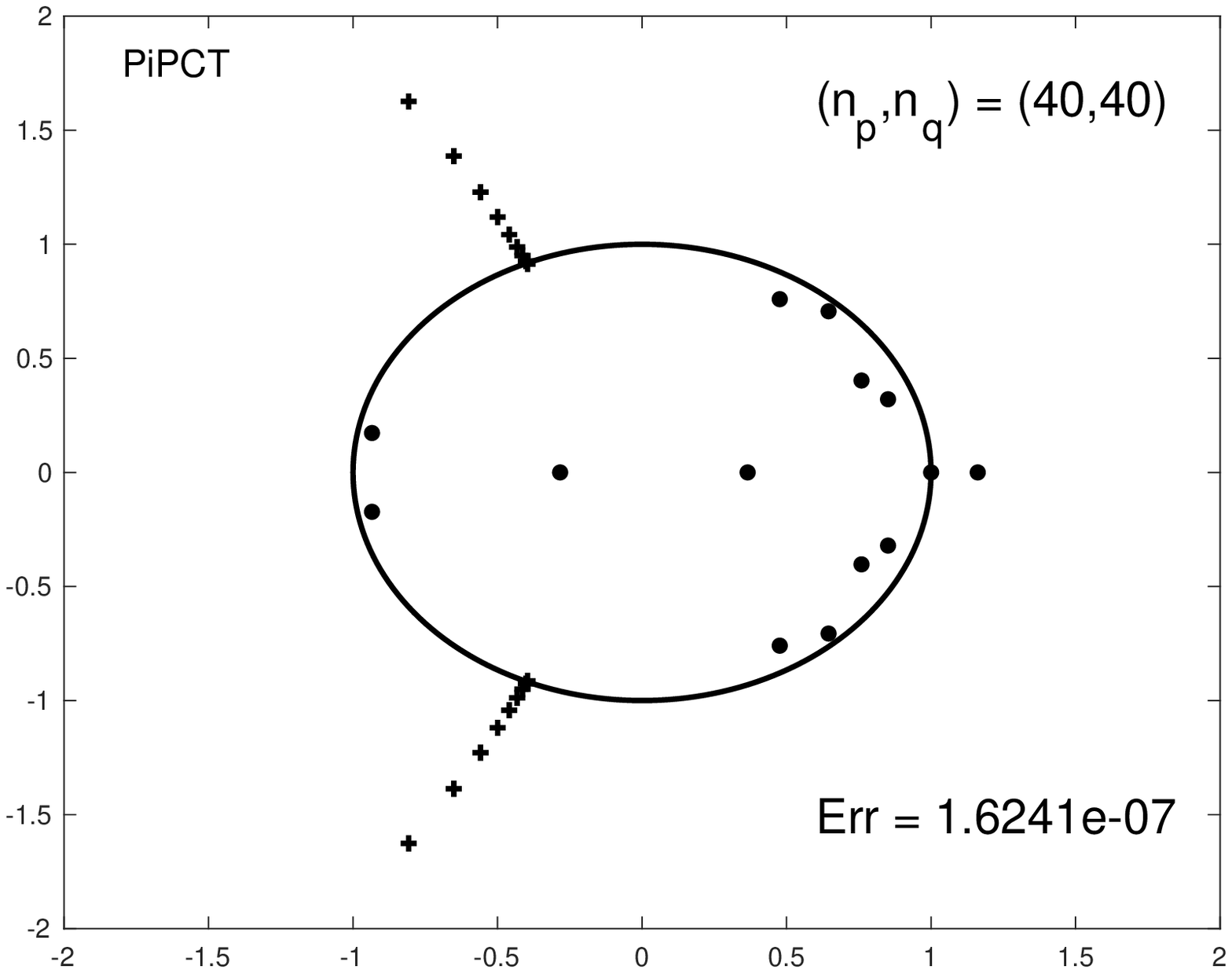}
	\includegraphics[height=5.2cm,width=5cm]{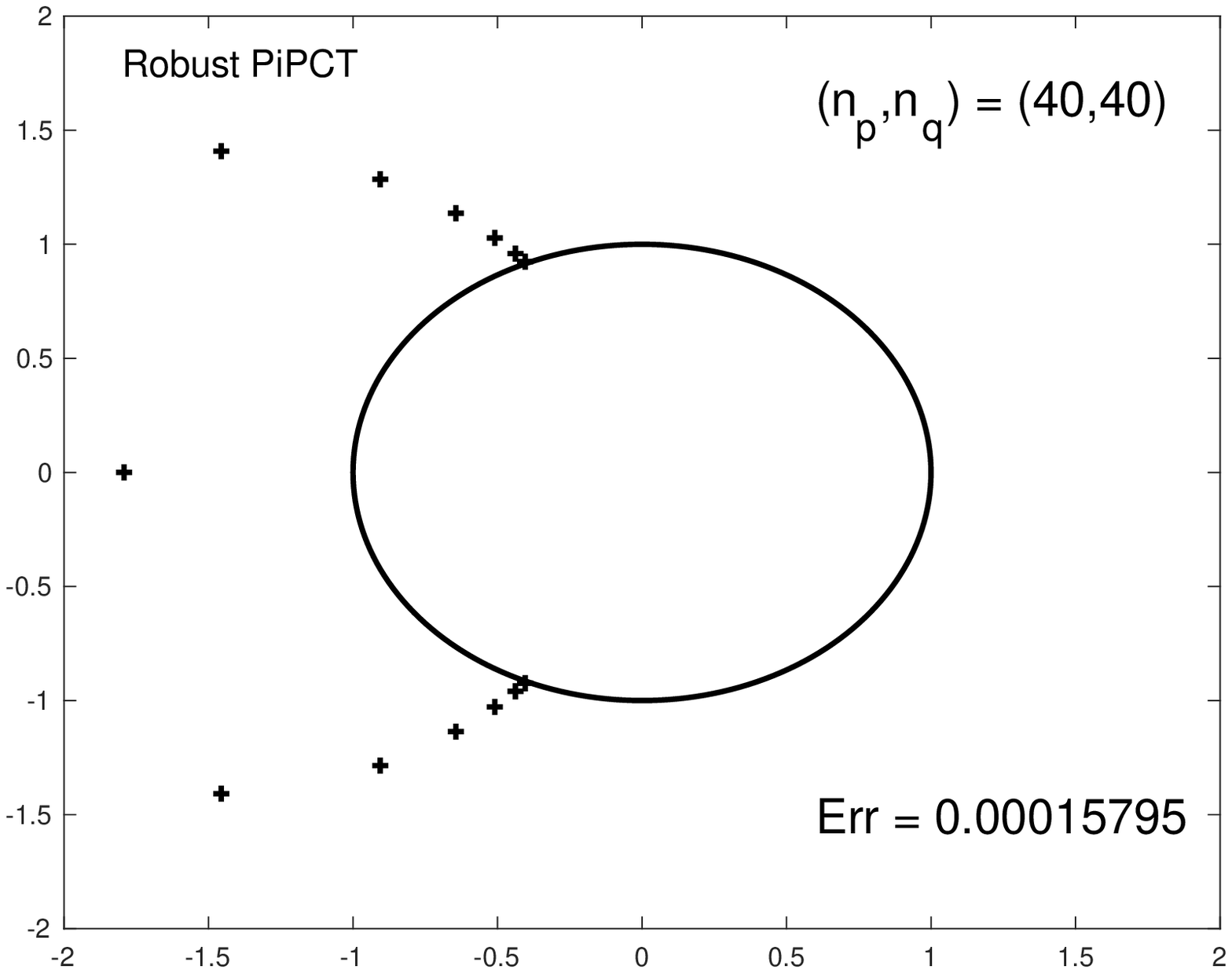}	
\caption{Poles of the PiPCT and Robust PiPCT (PiRPCT) methods in the vicinity of the jump discontinuity at $x=-0.4$ of $f$ given by  \eqref{discondeg.eq}.  Here, `+'-symbol denotes genuine poles and `$\bullet$'-symbol denotes spurious poles.}
\label{fig:RPCTpoles}\vspace{-0.1in}
\end{figure}

Figure \ref{fig:RPCTpoles} depicts the poles of PCT and RPCT approximations in $\epsilon$-badcell with jump discontinuity at $x = -0.4$. Plots in the first column of Figure \ref{fig:RPCTpoles} depict the poles of the PCT approximation, and the plots in the second column correspond to the poles of the RPCT approximation. Here `$\bullet$'-symbols denote the spurious poles, and `+'-symbols indicate genuine poles.  From the first column plots, we observe that the spurious poles do occur in the  $\epsilon$-badcells, which shows that the  APiPCT approximant is not free from Froissart doublets. We further observe that the number of spurious pole-zero pairs increases as $n_p(=n_q)$ increases, along with the increase in the accuracy.   Also, we observe from the second column plots of Figure \ref{fig:RPCTpoles} that the RPCT approximant is free from spurious poles, with a compromise in accuracy.  This is clearly because the RPCT method eliminates spurious poles (in this example) by reducing the degrees of the numerator and the denominator polynomials.

\section{Conclusion} $\!\!\!$
We proposed a piecewise Pad\'e-Chebyshev type (PiPCT) method for approximating piecewise smooth functions on a bounded interval. We obtained an $L^1$-error estimate in the regions where the target function is at least continuous. The error estimate shows that the order of accuracy depends on the degree of smoothness of the function as $N\rightarrow \infty$. We performed numerical experiments with a function involving both a jump discontinuity and a point at which the function is continuous but not differentiable. We demonstrated the convergence of the PiPCT method numerically in the vicinity of the singularities both in the case of increasing $N$ and increasing degrees of the numerator and the denominator polynomials. The numerical results clearly show the acceleration of the convergence of the approximant, which is not the case with the global Pad\'e-Chebyshev type approximation. The PiPCT method does not need the information about locations and types of singularities of the function. The PiPCT method is designed to work on a nonuniform mesh, which makes the algorithm more flexible for choosing a suitable adaptive partition and degrees. We used a strategy based on the zeros of the denominator polynomial to identify a cell, called a {\it badcell}, where a possible singularity of the function lies. Using the badcell identification, we further developed an adaptive way of partitioning the interval, which minimizes the computational time without compromising the accuracy. We also proposed a way to choose the numerator degree to improve the accuracy of the approximant in the bad cells, and the resulting algorithm is called the adaptive piecewise Pad\'e-Chebyshev type (APiPCT) algorithm. We performed numerical experiments to validate the APiPCT algorithm in terms of accuracy and time efficiency. We also compared the numerical results of PiPCT with some recently developed similar methods like the singular Pad\'e-Chebyshev method and the robust Pad\'e-Chebyshev method, and discussed the advantages of using PiPCT and APiPCT methods.
\section*{Acknowledgment}
The second author's work has been supported by the Science \& Engineering Research Board (SERB), Department of Science and Technology, under the project No. MTR/2017/000565.
\section*{References}


\end{document}